\documentclass[lettersize,journal]{IEEEtran}
\usepackage{amsmath,amsfonts}
\usepackage{algorithmic}
\usepackage{array}
\usepackage{textcomp}
\usepackage{stfloats}
\usepackage{url}
\usepackage{verbatim}
\usepackage{graphicx}
\def\BibTeX{{\rm B\kern-.05em{\sc i\kern-.025em b}\kern-.08em
    T\kern-.1667em\lower.7ex\hbox{E}\kern-.125emX}}
\usepackage{balance}

\usepackage{booktabs}
\usepackage{xcolor}
\usepackage{mathtools}
\usepackage{mathrsfs}

\usepackage{caption}
\usepackage{subcaption}

\usepackage{ulem}
\usepackage{float}

\newcommand{\bs}[1]{\boldsymbol{#1}}
\newcommand{\mb}{\mathbb}
\newcommand{\mc}{\mathcal}
\newcommand{\mr}{\mathrm}

\newcommand{\eps}{\epsilon}
\newcommand{\bsn}{{\boldsymbol{n}}} 
\newcommand{\bsx}{{\boldsymbol{x}}} 
\newcommand{\bsy}{{\boldsymbol{x'}}} 
\newcommand{\bsj}{{\boldsymbol{j}}} 
\newcommand{\bsm}{{\boldsymbol{m}}}

\newcommand{\bsE}{{\boldsymbol{E}}}
\newcommand{\bsH}{{\boldsymbol{H}}}

\newcommand{\mcE}{{\mathcal{E}}}

\newcommand{\bmk}{\bs{\mr{k}}}
\newcommand{\bg}{\bs{\gamma}}

\newtheorem{remark}{Remark}
\newenvironment{example}[1][]{\refstepcounter{example}\par\medskip
  \noindent \textbf{Example~\theexample. #1} \rmfamily}{\medskip}
\newcounter{example}[section]

\begin{document}

\title{Accelerated 3D Maxwell Integral Equation Solver using the Interpolated Factored Green Function Method}
\author{Jagabandhu Paul, \IEEEmembership{Member, IEEE} and Constantine Sideris, \IEEEmembership{Senior Member, IEEE}
\thanks{The authors gratefully acknowledge support by the Air Force Office of Scientific Research (FA9550-20-1-0087) and the National Science Foundation (CCF-2047433). J. Paul and C. Sideris are with the Department of Electrical and Computer Engineering, University of Southern California, Los Angeles, CA 90089, USA (e-mails: jagabandhu.paul@usc.edu, csideris@usc.edu).}}

\maketitle

\begin{abstract}
This article presents an \(\mc{O}(N\log N)\) method for numerical
solution of Maxwell's equations for dielectric scatterers using a
3D boundary integral equation (BIE) method. The underlying BIE method used is based on a
hybrid Nystr\"{o}m-collocation method using Chebyshev polynomials. It
is well known that such an approach produces a dense linear system,
which requires \(\mc{O}(N^2)\) operations in each step of an iterative solver. In this work, we propose an approach using the recently introduced Interpolated
Factored Green's Function (IFGF) acceleration strategy to reduce the cost of each iteration to \(\mc{O}(N\log N)\). To the best of our knowledge, this paper presents the first ever application of the IFGF method to fully-vectorial 3D Maxwell problems. The Chebyshev-based integral solver and IFGF method are first introduced, followed by the extension of the scalar IFGF to the vectorial Maxwell case. Several examples are presented verifying the \(\mc{O}(N\log N)\) computational complexity of the approach, including scattering from spheres, complex CAD models, and nanophotonic waveguiding devices. In one particular example with more than 6 million unknowns, the accelerated IFGF solver runs 42x faster than the unaccelerated method.
\end{abstract}



\section{Introduction\label{sec:intro}}
Solution of Maxwell's equations is quintessential to many modern
applications including antennas, microwave and photonic devices. The
lack of analytical solutions for general domains makes fast and 
accurate numerical solution methodology for Maxwell's equations
of utmost importance in such applications, in particular, for inverse
design approaches, where accurate field and gradient information is
needed in each iteration. The numerical methods available in the
literature can largely be classified into the following three groups:
finite-difference (FD) methods, finite-element (FE) methods,
and integral equation (IE) methods. 
For this work, we have used a boundary integral equation (BIE)
formulation due to its advantages for electromagnetic (EM) scattering problems over
FD~\cite{taflove_2005} and
FE~\cite{lalau_keraly_2013} methods, which are popular in the existing 
literature. In particular, BIE approaches only require discretizing
the surfaces of domains unlike the finite-difference and the
finite-element methods, which are volumetric in nature. In addition, BIE methods are almost dispersion-less due to analytically propagating information from sources to targets using Green's functions, unlike FD and FE methods. A BIE
approach was recently applied for
designing and optimizing photonic
devices~\cite{sideris_acs_2019,garza_acs_2023}, demonstrating significant improvements in both runtime and accuracy.  

The Method of Moments (MoM) approach is the most popular approach for discretizing BIEs in the
available literature. 
The authors in the pioneering work~\cite{rwg_1982} introduced the RWG basis
functions in order to solve the electric field integral equation in
conjunction with the MoM over flat triangular discretizations. A number of efforts have been made to alleviate some of the limitations arising from the first order basis
functions and improve the performance, including using higher order
basis functions (e.g.,~\cite{jorgensen_2004}) and phase-extracted basis functions (PEBFs)~\cite{pebf_dgie_2024}. Other high-order approaches also have been
introduced in the EM scattering context, such as~\cite{ganesh_2008}, which produces a spectrally accurate approximation of the tangential surface current using
a new set of tangential basis functions, as well as multiple other approaches as discussed in~\cite{notaros_2008}.

Recently, several works based on Nystr\"{o}m methods have been
proposed~\cite{rect_polar_bruno_garza_jcp_20}\cite{hu_cheby_tap_2021}\cite{garza_acs_2023}\cite{bruno_small_iter}\cite{notaros_2008}. In~\cite{rect_polar_bruno_garza_jcp_20},
the authors present a 
high-order method which decomposes the surface using non-overlapping
curvilinear parametric patches, after which each patch is discretized
using a Chebyshev grid in each parametric direction approximating the
unknown density using Chebyshev polynomials. This work was extended to
metallic and dielectric Maxwell scattering problems by leveraging the
Magnetic Field Integral Equation (MFIE) and the N-M\"{u}ller
formulations of the electromagnetic problems, respectively,
in~\cite{hu_cheby_tap_2021}. The method was further accelerated by
using GPU 
programming in~\cite{garza_acs_2023} for the indirect N-M\"{u}ller formulation and used to efficiently inverse design large 3D nanophotonic devices.
 
One of the main challenges of the Nystr\"{o}m method is that it produces dense linear systems leading to a \(\mc{O}(N^2)\)
computational complexity for iterative solvers once the integral operator has been discretized. In order to reduce this computational cost various algorithmic
acceleration strategies have been
proposed~\cite{fmm_3d_2006}~\cite{bruno_2001_equiv}. While these
methods are effective in reducing the asymptotic computational cost, they all rely on the Fast Fourier Transform (FFT), which presents
challenges for parallelization in the context of distributed memory parallel computer architectures. The recently introduced ``Interpolated Factored Green's
Function'' (IFGF) method~\cite{ifgf_jcp}, which relies on recursive
interpolation by means of Chebyshev expansion of relatively low degrees, on the other hand, does not use the FFT and therefore is immune to this
issue~\cite{ifgf_jcp_par}. Note that due to the low degree approximations used, the IFGF method does not yield spectral accuracy.

We leverage the high-order Nystr\"om method introduced in~\cite{rect_polar_bruno_garza_jcp_20,hu_cheby_tap_2021} to discretize the integral operators and repurpose the IFGF method (initially demonstrated for the scalar Helmholtz case in~\cite{ifgf_jcp}) to the 3D Maxwell scenario for accelerating the far interactions. To the best of our knowledge, this is
the very first application of the IFGF method to 3D Maxwell
electromagnetic scattering problems.

The rest of the paper is organized as follows: Section~\ref{sec:int_eqn_form} briefly reviews the indirect
N-M\"{u}ller formulation. Section~\ref{sec:numer_meth} presents
the numerical methodology used for the approximation of the integral
operators including the Chebyshev-based integral equation (CBIE)
method~\cite{rect_polar_bruno_garza_jcp_20,hu_cheby_tap_2021} used to decompose the
surface and for approximating the singular integrals, and 
the IFGF algorithm~\cite{ifgf_jcp} used to accelerate the far interactions. Section~\ref{sec:numer_meth} also presents the necessary details for extending the IFGF algorithm to the full-vectorial 3D Maxwell scenario, including the application of the IFGF algorithm to evaluate the integral operators for multiple densities and computation of the corresponding normal derivative of the single layer operator.
Finally, Section~\ref{sec:numer_exp} presents multiple numerical experiments
that demonstrate the performance and accuracy of the proposed method. In particular, we observe a \(\mc{O}(N\log N)\) complexity in
our calculations, and a significant increase in speed in approximation of the
operators. For instance, for a discretization containing more than 1.5 million points (more than 6 million unknowns), the computation time is reduced by a factor of \(42\) for one forward
map compared to the unaccelerated CBIE method. We conclude with a summary in Section~\ref{sec:summary}. 

\section{Integral Equation Formulation}\label{sec:int_eqn_form}
\begin{figure}
    \centering
    \includegraphics[scale=0.6]{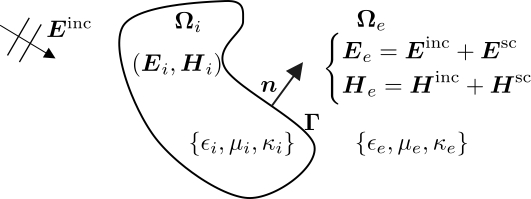}
    \caption{\footnotesize{Illustration for EM scattering by dielectric medium.}\label{fig:medium_illst}}
\end{figure}

We consider the numerical approximation of the scattered
electromagnetic fields, namely, the electric field \(\bsE\) and the
magnetic field \(\bsH\) in a dielectric media in three dimensions. For
simplicity, we consider two dielectric materials  
occupying the interior and exterior regions \(\Omega_{i}\) and
\(\Omega_{e}\), respectively, of the surface \(\Gamma\).
The case for a composition of domains can be treated
analogously~\cite{volakis_sertel_ie_em}. The
electromagnetic fields (\(\bsE, \bsH\)) and the material properties in
the interior 
and the exterior regions are denoted using the subscripts \(i\) and
\(e\), respectively, as illustrated in Fig.~\ref{fig:medium_illst}.
Suppose we are given the incident field
(\(\bsE^{\mr{inc}},\bsH^{\mr{inc}}\)), and that the exterior
unknown is the scattered field (\(\bsE^{\mr{sc}},\bsH^{\mr{sc}}\)),
where \(\bsE_e=\bsE^{\mr{inc}}+\bsE^{\mr{sc}}\) and
\(\bsH_e=\bsH^{\mr{inc}}+\bsH^{\mr{sc}}\). Then the equations 
\begin{equation}\label{eq:maxwell_int}
  \begin{cases}
    \nabla\times E_i- i\omega\mu_i H_i = 0,\\
    \nabla\times H_i + i\omega\epsilon_i E_i = 0,
  \end{cases}
\textrm{ in } \Omega_i;
\end{equation}
\begin{equation}\label{eq:maxwell_ext}
  \begin{cases}
    \nabla\times E_e- i\omega\mu_e H_e = 0,\\
    \nabla\times H_e + i\omega\epsilon_e E_e = 0,
  \end{cases}
\textrm{ in } \Omega_e;
\end{equation}
and the transmission 
conditions
\begin{equation}\label{eq:transmission}
  \begin{cases}
    \bsE_i\times \bsn = (\bsE^{\mr{sc}} + \bsE^{\mr{inc}})\times \bsn,\\
    \bsH_i\times \bsn = (\bsH^{\mr{sc}} + \bsH^{\mr{inc}})\times \bsn,
  \end{cases}
  \mathrm{ on }\; \Gamma,
\end{equation}
are satisfied along with the Silver-M\"{u}ller radiation condition.
An integral equation representation for the corresponding scattering problem can be found, for
instance, in~\cite[Chap. 5]{nedelec_aee},
\cite[Chap. 3]{volakis_sertel_ie_em} and~\cite{bruno_small_iter}.
As mentioned previously, we use a
three-dimensional boundary integral
representation to obtain the
electromagnetic fields~\cite{colton_kress_inverse_em}. In particular, we consider the indirect
N-M\"{u}ller formulation~\cite{c_muller}. To define the integral
representation, let us define the tangential integrals
\begin{align}
  \bs{R}_{s}[\bs{\Phi}](\bsx) &= \bsn\times\int_{\Gamma}
                                G_s(\bsx,\bsy) \bs{\Phi}(\bsy)
                                d\gamma(\bsy),
                                \label{eq:potential_r}\\
  \bs{K}_{s}[\bs{\Phi}](\bsx) &= \bsn\times\nabla\times\int_{\Gamma} G_s(\bsx,\bsy) \bs{\Phi}(\bsy)
               d\gamma(\bsy),\label{eq:potential_k}\\
  \bs{T}_{s}[\bs{\Phi}](\bsx) &= \bsn\times\int_{\Gamma}\nabla G_s(\bsx,\bsy) \mathrm{div}_{\Gamma}\bs{\Phi}(\bsy) d\gamma(\bsy)\label{eq:potential_t}
\end{align}
where \(s\) denotes the domain subscripts \(e\) or \(i\), and  \(G_{s}(\bsx,\bsy)=\exp(i\kappa_{s}|\bsx-\bsy|)/(4\pi|\bsx-\bsy|)\)
is the homogeneous space Green's function of the Helmholtz equation. The
indirect N-M\"{u}ller integral formulation can then be written in matrix from as~\cite{garza_acs_2023}
\begin{equation}\label{eq:n_muller}
\begin{split}
&\begin{bmatrix}
  \frac{(\mu_e+\mu_i)}{2}\bs{I} + \bs{K}^{\Delta}_{\mu} & -(\bs{R}^{\Delta}
                                                       + \bs{T}^{\Delta})\\
  (\bs{R}^{\Delta}+ \bs{T}^{\Delta}) & \frac{(\eps_e+\eps_i)}{2}\bs{I} + \bs{K}^{\Delta}_{\eps}
\end{bmatrix}
\begin{bmatrix}
  \bsm \\
  \bsj
\end{bmatrix}  \\
&\quad\quad =
\begin{bmatrix}
  \omega^{-1}\, \bsn\times \bsE^{\mr{inc}} \\
  \omega^{-1}\, \bsn\times \bsH^{\mr{inc}}
\end{bmatrix}
  \end{split}
\end{equation}
where
\begin{align}\label{eq:k_delta}
   \bs{K}^{\Delta}_{\alpha} &\equiv (\alpha_{e}\bs{K}_{e}-\alpha_{i}\bs{K}_{i}),\\
   \bs{R}^{\Delta} &\equiv
-i\omega(\mu_{e}\epsilon_{e}\bs{R}_{e} - \mu_{i}\epsilon_{i}\bs{R}_{i})
\equiv \frac{-i}{\omega}(\kappa^{2}_{e}\bs{R}_{e}-\kappa^{2}_{i}\bs{R}_{i}),\label{eq:r_delta}\\
   \bs{T}^{\Delta} &\equiv \frac{-i}{\omega}(\bs{T}_{e}-\bs{T}_{i}),\label{eq:t_delta}
\end{align}
with \(\bs{I}\) denoting the identity operator, and \(\alpha\)
denoting either \(\epsilon\) and \(\mu\) (see Remark~\ref{rmk:weakly_sing}). 
\begin{remark}[\textbf{Weak singularity in~\eqref{eq:n_muller}}]\label{rmk:weakly_sing}
\normalfont{The first two terms in the expansion of the exponential
function in (\(\nabla G_{e}-\nabla G_{i}\))
vanish~\cite{oijala_muller_tap}, thus the operators in~\eqref{eq:n_muller} 
remain weakly singular.}
\end{remark}

\section{Numerical Methodology\label{sec:numer_meth}}
As discussed in~\cite{garza_thesis}, the integral operators in~\eqref{eq:n_muller}
can be re-expressed via algebraic manipulation in terms of the weakly-singular
operators of the form
\begin{equation}\label{eq:sl_scalar}
    S[\phi](\bsx) = \int_{\Gamma} G_{\kappa}(\bsx,\bsy) \phi(y) d\gamma(\bsy),
\end{equation}
and the normal derivative
\begin{equation}\label{eq:dl_scalar}
    D[\phi](\bsx) = \frac{\partial S[\phi](\bsx)}{\partial n(\bsx)},
\end{equation}
where \(\phi\) represents a Cartesian component of one of the surface current
densities \(\bsj=(j_{x},j_{y},j_{z})\) and
\(\bsm=(m_{x},m_{y},m_{z})\); or one of the surface divergences
\(\mr{div}_{\Gamma}\bsj\) and \(\mr{div}_{\Gamma}\bsm\) of the current
densities. Note that for the surface divergences
\(\mr{div}_{\Gamma}\bsj\) and \(\mr{div}_{\Gamma}\bsm\), only the
operator \(S[\phi]\) needs to be computed. The EM 
operators are then evaluated by further compositions of other 
operators such as gradient as per the
definitions of the operators~\eqref{eq:k_delta}-\eqref{eq:t_delta} involved in~\eqref{eq:n_muller}.
The required derivatives are obtained 
using the Chebyshev recurrence relations~\cite{press_num_recipes} based on the underlying
Chebyshev grids given by~\eqref{eq:cheby_nodes} below. To approximate the above integrals and to
further obtain accurate values of the EM operators, the proposed
IFGF accelerated Chebyshev-based (CBIE-IFGF) method
discretizes the surface \(\Gamma\) as described in what follows.

\subsection{CBIE-IFGF: Surface Discretization.}\label{sec:surf_discr}
\begin{figure}
\includegraphics[scale=0.5]{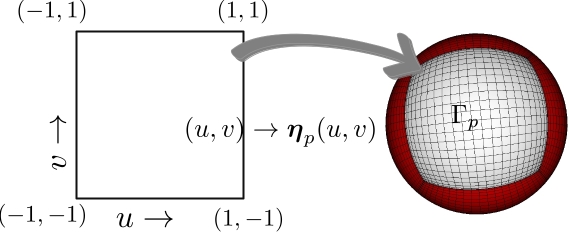}
\caption{\footnotesize{Mapping the square \([-1,1]\times [-1,1]\) in parametric \((u,v)\)-space to a
patch on the surface of a sphere in Cartesian coordinates.}\label{fig:patch_cheby}}
\end{figure}

The proposed CBIE-IFGF method decomposes the surface \(\Gamma\) into a set of
\(P\) non-overlapping, curvilinear quadrilateral patches \(\Gamma_{p};
p=0,\dots,(P-1)\) such that
\begin{equation}\label{eq:surf_patch}
    \Gamma = \bigcup_{p=0}^{P-1} \Gamma_{p},
\end{equation}
where each patch \(\Gamma_{p}\) has a parametrization
\(\bs\eta_{p}(u,v)\) with \((u,v)\in [-1,1]^{2}\), as illustrated Fig.~\ref{fig:patch_cheby}.
Following~\eqref{eq:sl_scalar} and~\eqref{eq:surf_patch}, we can write
\begin{equation}\label{eq:sl_patch}
    S[\phi](\bsx) = \sum_{p=0}^{P-1}S_{p}(\bsx),
\end{equation}
where
\begin{equation}\label{eq:sl_p}
    S_{p}(\bsx) = \!\int\limits_{[-1,1]^{2}}\!\!\! G_{\kappa}(\bsx,\bs\eta_{p}(u,v))
\phi_{p}(u,v) J_{p}(u,v) dudv,
\end{equation}
\(G_{p}(\bsx,u,v)=G(\bsx,\bs\eta_{p}(u,v))\), \(\phi_{p}(u,v)=\)
\(\phi(\bs\eta_{p}(u,v))\), and 
\(J_{p}(u,v)\) denotes the Jacobian of the parametrization
\(\bs\eta_{p}\).

In each of the patches, the method places an open Chebyshev grid
along \(u\) and \(v\) directions:
\begin{equation}\label{eq:cheby_nodes}
    u_{i} = \cos((k+0.5)\pi / N_{C}): k=0,\dots,N_{C}-1,
\end{equation}
(similarly, for \(v\)) for certain integer value of \(N_{C}\).
(For simplicity of the description, we have used
the same number of points in both 
parametric directions---which is not a necessity.) The collection \(\Gamma_{N}\)
of the discretization points on \(\Gamma\) is given below:
\begin{align}\label{eq:mesh_pts_l}
  \Gamma_{N} = 
  \left\{\right.
    \bsx_{\ell}=&\left. \bs\eta_{p}(u_{i},v_{j}):
    \ell=j+N_{C}i+pN^{2}_{C}; \right.\nonumber \\
       &\left. 0\leq i,j \leq N_{C}-1; 0\leq p
         \leq P-1 \right\}.
\end{align}
The total number of discretization points are \(N=P\times
N_{C}^{2}\). 
The resulting linear system that arises upon discretizing~\eqref{eq:n_muller} 
based on the discretization \(\Gamma_{N}\) of the surface \(\Gamma\)
is denoted by
\begin{equation}\label{eq:lin_system}
    A\bs{y} = b,
\end{equation}    
where the matrix \(A\) is of size \( 4N\times 4N\). The proposed method uses
GMRES to solve the linear system in~\eqref{eq:lin_system}. Next, we
discuss the approximation methodology for the integral \(S_{p}(\bsx)\)
over a patch \(\Gamma_{p}\) for \(p=0,\dots,P-1\) for a given target point \(\bsx\).

\subsection{Approximating the Integral \(S_{p}(\bsx)\) for \(p=0,\dots,P-1\).}
There are two main difficulties in approximating
\(S_{p}(\bsx)\)~\eqref{eq:sl_p} for \(p=0,\dots,P-1\):
first, the singular and near-singular behavior of the kernel when
the target point \(\bsx\) and the source point \(\bsy=\bs\eta_{p}(u,v)\) coincide or,
are in close proximity; and second, the \(\mc{O}(N^{2})\) cost required
to evaluate the discrete integrator for the regular/non-singular case when the
source and the target points are sufficiently away from each other.
In order to classify the singular and non-singular target/observation
points for integration over a given surface patch \(\Gamma_p\), we
define the following 
``distance of a point \(\bsx\) from a surface patch \(\Gamma_{p}\)'':
\begin{equation}\label{eq:dist}
    \mathrm{dist}(\bsx,\Gamma_{p}) = \inf\{|\bsx-\bsy| : \bsy\in \Gamma_{p}\}.
\end{equation}
Moreover, for a certain number \(\delta>0\), let us define the index sets
\(\mathscr{S}_{\ell}\) and \(\mathscr{R}_{\ell}\) for a point \(\bsx_{\ell}\)
\begin{align}
  \mathscr{S}_{\ell} &\coloneq \{p:
                      \mathrm{dist}(\bsx_{\ell},\Gamma_{p})\leq \delta\},~\mathrm{and}\label{eq:patch_sing}\\
  \mathscr{R}_{\ell} &\coloneq \{p: \mathrm{dist}(\bsx_{\ell},\Gamma_{p}) > \delta\}\label{eq:patch_reg}
\end{align}
containing the indices of the patches \(\Gamma_{p}\) \((p=0,\dots,P-1)\) for which the point
\(\bsx_{\ell}\) is treated as a singular, or near-singular and a
regular/non-singular point, respectively.
Thus
\begin{equation}\label{eq:sl_patch_1}
    S[\phi](\bsx_{\ell}) = \sum_{p\in\mathscr{S}_{\ell}} S_{p}(\bsx_{\ell}) + \sum_{p\in\mathscr{R}_{\ell}} S_{p}(\bsx_{\ell}).
\end{equation}  
We start our discussion on the quadrature to approximate
\(S_{p}(\bsx_{\ell})\) 
 for the non-singular case  
(\(p\in\mathscr{R}_{\ell}\)) and the singular case (\(p\in\mathscr{S}_{\ell}\)) in 
Section~\ref{sec:non_sing} and Section~\ref{sec:rect_pol},
respectively, in what follows.

\subsection{Approximation of \(S_{p}(\bsx_{\ell})\): Non-Singular Case.\label{sec:non_sing}}
The non-singular case can be treated using Fej\'{e}r quadrature based on the nodes given
by~\eqref{eq:cheby_nodes} and the weights
\begin{equation}\label{eq:fejer_weighs}
    w_{n} = \frac{2}{N_{C}} \left[1 - 2
\sum\limits_{m=1}^{\lfloor{N_{C}/2}\rfloor} \frac{1}{4m^{2}-1}\cos\left(\frac{m\pi(2n+1)}{N_{C}}\right)\right],
\end{equation}    
for \(n=0,\dots,N_{C}-1\). A straightforward application of
the Fej\'{e}r quadrature at \(\bsx_{\ell}\) reads 
\begin{equation}\label{eq:non_sing}
    \mc{S}_{p}(\bsx_{\ell}) \approx
\sum\limits_{i=0}^{N_{C}-1}\sum\limits_{j=0}^{N_{C}-1}
G(\bsx_{\ell},u_{i},v_{j}) \phi_{p}(u_{i},v_{j}) J_p(u,v) w_{i}w_{j}.
\end{equation}
Clearly, such an application of the Fej\'{e}r quadrature has an
overall complexity of \(\mc{O}(N^{2})\); in order to get a 
reduced \(\mc{O}(N\log N)\) run time in the non-singular calculation,
it is accelerated using the IFGF algorithm discussed in
Section~\ref{sec:ifgf_acc}.


\subsection{Approximation of \(S_{p}(\bsx_{\ell})\): Singular Case.}\label{sec:rect_pol}
The evaluation of \(S_{p}(\bsx_{\ell})\) for the singular case requires
special treatment via change of variable to resolve the singularity
present in the kernel. Given that this change of variable as such
depends on the target point, in order to make the singular calculation
more efficient, the method uses precomputed integral moments
of the product of the kernel and the Chebyshev functions of the
parametric variables \(u\) and \(v\), along with the Jacobian of the
parametric change of variable. Letting the Chebyshev expansion of the
density \(\phi_{p}(u,v)\):
\begin{equation}\label{eq:phi_cheby}
   \phi_{p}(u,v) \approx \displaystyle\sum\limits_{m=0}^{N_{C}-1}
   \sum_{n=0}^{N_{C}-1}\!\! a^{p}_{m,n}T_{m}(u)T_{n}(v)
\end{equation}
based on the nodes~\eqref{eq:cheby_nodes}; using~\eqref{eq:sl_p} the
integral \(S_{p}(\bsx)\) can be approximated by evaluating the summation
\begin{equation}\label{eq:rp_exp}
  S_{p}(\bsx_{\ell})\nonumber \approx
   \displaystyle\sum\limits_{m=0}^{N_{C}-1}\sum_{n=0}^{N_{C}-1}
    a^{p}_{m,n} \beta^{p,\ell}_{m,n},
\end{equation}
where the moments have been computed as
\begin{equation}\label{eq:int_moments}
     \beta^{p,\ell}_{m,n} = \displaystyle\int\limits_{[-1,1]^{2}}
     G(\bsx_{\ell},u,v) T_{m}(u)T_{n}(v) J_{p}(u,v) dudv.
\end{equation}

To evaluate the moments \eqref{eq:int_moments} at a target point
\(\bsx_{\ell}\in\Gamma_{N}\), the method finds the point \(\bsx_{\ell}^{p}=\bsx_{p}(u^{p}_{\ell},v^{p}_{\ell})\in\Gamma_{p}\)
that is closest to \(\bsx_{\ell}\), or rather its parametric-space coordinates,
which in general can be found as the solution of the distance
minimization problem
\begin{equation}\label{eq:dist_min}
   (\tilde{u}^{p}_{\ell},\tilde{v}^{p}_{\ell}) =
   \underset{(u,v)\in [-1,1]^{2}}{\mathrm{arg\,min}}
   |\bsx_{\ell}-\bs{\eta}_{p}(u,v)|.
\end{equation}
Following \cite{rect_polar_bruno_garza_jcp_20}, the method utilizes
the golden section search algorithm, with initial bounds obtained from
a direct minimization over all of the original discretization points
in \(\Gamma_{p}\). Note that if \(\bsx_{\ell}\) is itself
a grid point in \(\Gamma_{p}\) and
\(\bsx_{\ell}=\bs\eta_{p}(u_{m},v_{n})\) then the parametric
coordinates \((\tilde{u}^{p}_{\ell},\tilde{v}^{p}_{\ell}) =
(u_{m},v_{n})\).

The method utilizes a one-dimensional change of variable
to each coordinate in the \(uv\)-space to construct a clustered grid
around a given target point. To this end, we consider the following
mapping \(w:[0,2\pi]\rightarrow [0,2\pi]\), with parameter \(d\geq 2\)
(see,~\cite[Section 3.5]{colton_kress_inverse_em}):
\begin{equation}\label{eq:mk_cov_1}
   w(\tau) = 2\pi
\frac{[\nu(\tau)^{d}]}{[v(\tau)]^{d}+[\nu(2\pi-\tau)]^{d}},\;
0\leq\tau\leq 2\pi,
\end{equation}
where
\begin{equation}\label{eq:mk_cov_2}
  v(\tau) = \left(\frac{1}{d}-\frac{1}{2}\right)
\left(\frac{\pi-\tau}{\pi}\right)^{3} +
\frac{1}{d}\left(\frac{\tau-\pi}{\pi}\right)+\frac{1}{2}. 
\end{equation}
It can be shown that \(w\) has vanishing derivatives up to order
\(d-1\) at the interval endpoints.

Further, the one-dimensional change of variables (defined on the basis
of the change of variable above)
\begin{equation}\label{eq:cov}
   \xi_{\alpha}(\tau) = 
  \begin{cases}
    \alpha + \left(\frac{\mathrm{sgn}(\tau)-\alpha}{\pi}\right) w(\pi
    |\tau|), \; \mathrm{ for }\; \alpha \neq \pm 1,\\
    \alpha - \left(\frac{1 + \alpha}{\pi}\right) w(\pi
    |\frac{\tau-1}{2}|), \; \mathrm{ for }\; \alpha = 1,\\
    \alpha + \left(\frac{1 - \alpha}{\pi}\right) w(\pi
    |\frac{\tau+1}{2}|), \; \mathrm{ for }\; \alpha = -1
  \end{cases}
\end{equation}    
clusters the points around \(\alpha\). Indeed, a use of Fej\'{e}r's
rule with \(N_\beta\) points yields 
\begin{eqnarray}\label{eq:int_moments_cov}
   \beta^{p,\ell}_{m,n}\approx
  \displaystyle\sum\limits_{i=0}^{N_{\beta}-1}\sum\limits_{j=0}^{N_{\beta}-1}\!\!
  &\!\!G(\bsx_{\ell},u_{i}^{p,\ell},v_{j}^{p,\ell})
    J_{p}(u_{i}^{p,\ell},v_{j}^{p,\ell}) T_{m}(u_{i}^{p,\ell})\nonumber \\
  &T_{n}(v_{j}^{p,\ell}) \mu_{i}^{u,p,\ell}\mu_{j}^{v,p,\ell} w_{i}w_{j}
\end{eqnarray}
where 
\begin{eqnarray}\label{eq:int_moments_1}
  u_{i}^{p,\ell} = \xi_{\tilde{u}_{\ell}^{p}}(t_{i}),&
  v_{j}^{p,\ell}=\xi_{\tilde{v}_{\ell}^{p}}(t_{j}),\\
  \mu_{i}^{u,p,\ell} =
  \frac{d\xi_{\tilde{u}_{\ell}^{p}}}{d\tau}(t_{i}),&
  \mu_{j}^{v,p,\ell} =
    \frac{d\xi_{\tilde{v}_{\ell}^{p}}}{d\tau}(t_{j}),
\end{eqnarray}
for \(i,j=0,\dots,N^{\beta}-1\). These integral moments (and the
required closest points) are computed only once at the beginning of each run.
The value of \(N_{\beta}\) is chosen (based on the experiments
presented in~\cite[Fig.~3a-3b and TABLE II]{hu_cheby_tap_2021}) to match the accuracy provided
by the IFGF interpolation algorithm, which we describe in the following subsection.

\subsection{IFGF Acceleration of Non-singular Calculation.}\label{sec:ifgf_acc}
As discussed in Section~\ref{sec:non_sing},
for a given surface discretization \(\Gamma_N:=\{x_1 ,\dots, x_N\}\),
we need to evaluate a discrete sum of the form 
\begin{equation}\label{eq:ifgf_sum}
  I(x_\ell) := \sum_{\substack{m=1 \\ m\neq\ell}}^{N} a_m G(\bsx_\ell,\bsx_m), \quad \ell=1,\dots,N
\end{equation}
where \(N\) denotes a given positive integer, \(m=1,\dots,N, \) and,
 \(\bsx_m \in\Gamma_{N}\) and \(a_m \in \mb{C}\) denote pairwise
different points and given complex numbers, respectively
(cf.~\eqref{eq:sl_patch_1} for \(p\in\mathscr{R}_{\ell}\)
and~\eqref{eq:non_sing}). Clearly, a direct
evaluation of the sum \(I(x)\) for all \(x\in\Gamma_N\) requires
\(\mathcal{O}(N^2)\) operations. The recursive interpolation based
IFGF approach, introduced in~\cite{ifgf_jcp}, can reduce the
cost of this evaluation to \(\mathcal{O}(N\log N)\).
We first describe the IFGF interpolation strategy for a scalar density
\(\phi\), 
 and for the
contribution coming from the source points within a certain box in what follows.


\subsection{IFGF Interpolation: Single Source Box.}\label{sec:ifgf_single_box}
Consider the axis aligned box \(B(\bsx^{S},H)\) 
\begin{align*}
 B&\left(\bsx^{S},H\right):=
    \left[(\bsx_{S})_1-\frac{H}{2},(\bsx_{S})_1+\frac{H}{2}\right)\times\\
    &\left[(\bsx_{S})_2-\frac{H}{2},(\bsx_{S})_2+\frac{H}{2}\right)\times
      \left[(\bsx_{S})_3-\frac{H}{2},(\bsx_{S})_3+\frac{H}{2}\right)
\end{align*}
in \(\mathbb{R}^3\) of side length \(H>0\) and centered at
\(x^{S}=((x_{S})_1,(x_{S})_2,(x_{S})_3)\in\mathbb{R}^3\). 
Let
\( \Gamma^S_{N_S} := \{x_1^S, \dots, x_{N_S}^S\} \in B(x^S, H) \cap
 \Gamma_N \) denote the source points contained in \(B(x^{S},
H)\). Then, the contribution of the source points within the box
\(B(x^S, H)\) at a target point \(\bsx_{\ell}^{T}\) is given by
\begin{equation}\label{eq:box_contr}
    I_S(\bsx_\ell^T) := \displaystyle\sum_{m=1}^{N_S} a_m^S G(\bsx_\ell^T, \bsx_m^S), \quad \ell=1,\dots,N_T.
\end{equation}
To efficiently evaluate the sum \(I_{S}(\bsx)\), the IFGF
method first factorizes the kernel as 
\begin{equation}\label{eq:ifgf_fact}
G(\bsx,\bsy) = G(\bsx,\bsx^{S}) g_{S}(\bsx,\bsy)
\end{equation}
where
\(G(\bsx,\bsx^{S}) =
\frac{e^{i\kappa|\bsx-\bsx^{S}|}}{4\pi|\bsx-\bsx^{S}|}\) is called the
``centered factor'', and 
\begin{align}\label{eq:ifgf_ana_factor}
  g_S(\bsx,\bsy) &= \frac{|\bsx-\bsx^{S}|}{|\bsx-\bsy|} e^{i\kappa(|\bsx-\bsy|-|\bsx-\bsx^{S}|)}
\end{align}
is called the ``analytic factor''. Using the
factorization~\eqref{eq:ifgf_fact}, we can write~\eqref{eq:box_contr} as
\begin{equation}\label{eq:box_contr_factor}
    I_S(\bsx_{\ell}^{T}) =
  \displaystyle\sum_{m=1}^{N_S} a_m^S G(\bsx_{\ell}^{T},\bsx_m^S) = G(\bsx_{\ell}^{T},\bsx_S) F_S(\bsx_{\ell}^{T})
\end{equation}
where
\begin{equation}\label{eq:box_analfac}
    F_S(\bsx_{\ell}^{T}) = \displaystyle\sum_{m=1}^{N_S} a_m^S g_S(\bsx_{\ell}^{T},\bsx_m^S)
\end{equation}
is the analytic factor contribution.

\begin{figure}[ht]\label{fig:kernel_plot}
\centering
\includegraphics[scale=0.5]{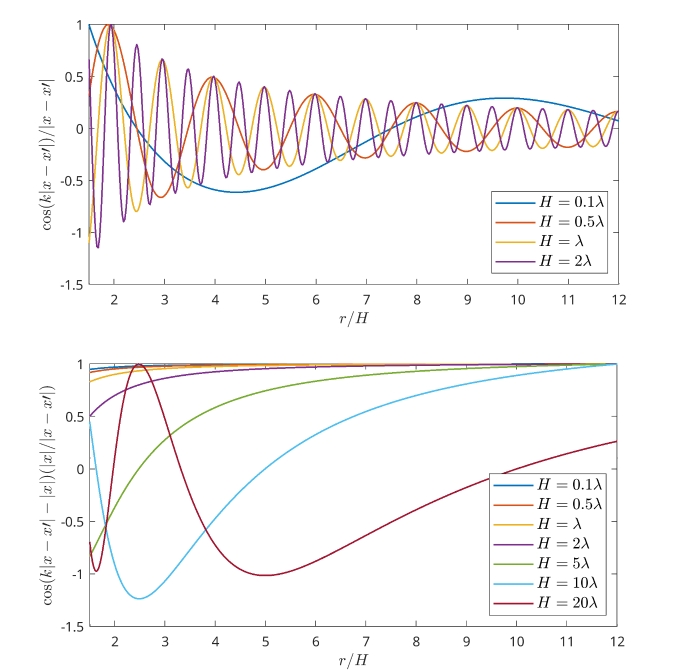}
\caption{\footnotesize{Source point factorization test. Consider an
origin centered box of side \(H\). In this illustration, the source 
point is positioned at \(\bsx'=(0,0,H/2)\), 
and the point \(\bsx\) varies on the line \(t(1,0,0)\) for
\(t\in[3H/2,12H]\). The graph at the top displays the real part
of the kernel \(G(\bsx,\bsy)\), which illustrates a highly-oscillatory behavior even
for the smaller size boxes. The bottom graph illustrates the slow
oscillatory behavior of the analytic factor \(g_{S}(\bsx,\bsy)\) for boxes as large as
\(H=20\lambda\).}\label{fig:kernel_plots}} 
\end{figure}

Defining the radius of the source box \(B(\bsx^{S},H)\) by
\(h:=\frac{\sqrt{3}}{2}H\), and the variable \(s:=h/r\) where \(r:=|\bsx-\bsx^{S}|\),
the IFGF method then uses the change of variables
\begin{equation}\label{eq:ifgf_cov}
    \bsx = x(s,\theta,\phi;\bsx^{S}) = \widetilde{x}(h/s,\theta,\phi;\bsx^{S}),
\end{equation}
where \(\widetilde{x}\) denotes the \(\bsx^{S}\)-centered spherical
coordinate change of variables
\begin{equation}
  \widetilde{x}(r,\theta,\phi;\bsx^{S}):= \bsx^{S} + \begin{pmatrix}
    r\sin\theta\cos\phi \\
    r\sin\theta\sin\phi \\
    r\cos\theta
  \end{pmatrix};
  0 \leq r=|\bsx-\bsx^{S}| < \infty
\end{equation}
with
\(0 \leq \theta \leq \pi, 0 \leq \phi < 2\pi\).
Using the variables defined above, we can write
\(G(\bsx,\bsx^{S})=\exp(i\kappa r)/(4\pi r)\), and 
\begin{equation}\label{eq:ana_factor_sphr}
    g_S(\bsx,\bsy) =
 \frac{\exp\left(i\kappa r\left(\left|\frac{\bsx-\bsx^{S}}{r}-\frac{\bsy-\bsx^{S}}{h}s\right|-1\right)\right)} {\left|\frac{\bsx-\bsx^{S}}{r} - \frac{\bsy-\bsx^{S}}{h}s\right|}.
\end{equation}
Since by definition \(|\bsx-\bsx^{S}|=r\),
and for any point \(\bsy\in B(\bsx^{S},H)\)---\(|\bsy-\bsx^{S}|\leq h\), the 
function \(g_{S}\) is analytic for all \(|s| < 1\), including \(s=0\)
(\(r=\infty\))~\cite{ifgf_jcp}. 
Moreover, the function \(g_{S}\) is slowly
oscillatory in the \(s\) variable, as shown in Fig.~\ref{fig:kernel_plots}.


As discussed in the next section, the IFGF method
approximates the contribution \(I_{S}(\bsx_{\ell}^{T})\)
only at target points which are at least one box
away. For such target points, we have 
\((s,\theta,\phi)\in \mc{E} = [0,\eta]\times [0,\pi] \times [0,2\pi)\), where
\(\eta=\sqrt{3}/3< 1\). Hence, using the analytic and low-oscillatory
characteristics of \(g_{S}\), and by linearity, the
analytic factor contribution \(F_{S}(\bsx_{\ell}^{T})\) can
efficiently be interpolated in the \((s,\phi,\theta)\) variables with
a few finite interpolation intervals in the \(s\) variable, in
conjunction with a number of adequately selected compact interpolation
intervals in the \(\phi\) and \(\theta\) variables. Once
\(F_{S}(\bsx_{\ell}^{T})\) has been interpolated, the value of
\(I_{S}(\bsx_{\ell}^{T})\) can be obtained by multiplying the centered
factor \(G(\bsx_{\ell}^{T},\bsx^{S})\), which is known analytically.

\begin{figure}
  \centering
\includegraphics[scale=0.5,trim={20 20 20 20},clip]{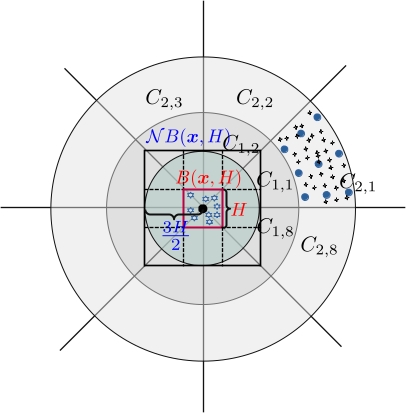}
\caption{\footnotesize{Illustration of a cone structure for a single box in 2D. The figure
illustrates a cone-structure for the box
\(B(\bsx,H)\) shown in red, the neighborhood boxes in \(\mc{N}B(\bsx,H)\) separated by dashed lines 
within thick black lines, and a \(P_s\times P_a=3\times 3\) interpolation
scheme in cone-segment \(C_{2,1}\). The hollow blue star-shape points
denote the source points within \(B(\bsx,H)\), blue solid circles
denote the interpolation grid points in \(C_{2,1}\), and the points
denotes by black-\(+\) signs are the target points within \(C_{2,1}\). The
IFGF method does not interpolate contribution from \(B(\bsx,H)\) at
any target point that  
lies within \(\mc{N}B(\bsx,H)\).}\label{fig:cone_struct_2d}}
\end{figure}

To appropriately structure the interpolation, the domain \(\mcE\) is
partitioned as described in what follows. For given positive integers
\(n_{s}\) and \(n_{a}\), we denote the length of the interpolation
intervals in \(s\), and in the angular variables \(\theta\)
and \(\phi\) as
\begin{equation}
  \Delta_{s}=\frac{\eta}{n_{s}}\;\mathrm{ and }\; \Delta_{\theta}=\Delta_{\phi}=\frac{\pi}{n_{a}},
\end{equation}
respectively. Then for each \(\gamma=(\gamma_{1},\gamma_{2},\gamma_{3})\) in 
\(K_{C}:= \{1,\dots,n_{s}\} \times \{1,\dots,n_{a}\}\times \{1,\dots,2 n_{a}\},\)
the interpolation intervals
\(E^{s}_{\gamma_{1}}\), \(E^{\theta}_{\gamma_{2},\gamma_{3}}\)
and \(E^{\phi}_{\gamma_{3}}\) along the \(s\), \(\theta\),
and \(\phi\) directions are defined by 
\begin{equation}
    E^{s}_{\gamma_{1}} = [(s-1)\Delta_{s},s\Delta_{s})
\end{equation}
\begin{equation}
    E^{\theta}_{\gamma_{2},\gamma_{3}} = 
\begin{cases}
  [(n_{a}-1)\Delta_{\theta},\pi]
  &\mathrm{for}\;
    \gamma_{2}=n_{a},\gamma_{3}=2n_{a}\\
  (0,\Delta_{\theta})
  &\mathrm{for}\;
    \gamma_{2}=1,\gamma_{3}>1\\
  [(\gamma_{2}-1)\Delta_{\theta},\gamma_{2}\Delta_{\theta}]
  &\mathrm{otherwise},
\end{cases}
\end{equation}
\begin{equation}
    E^{\phi}_{\gamma_{3}} = [(\gamma_{3}-1)\Delta_{\phi},\gamma_{3}\Delta_{\phi}).
\end{equation}
For each  \(\gamma=(\gamma_{1},\gamma_{2},\gamma_{3})\in K_{C}\), we
call
\begin{equation}\label{eq:cone_dom}
E_{\bs\gamma} := E^{s}_{\gamma_{1}}\times E^{\theta}_{\gamma_{2}}\times E^{\phi}_{\gamma_{3}}
\end{equation}
a \emph{cone domain}, and the image of \(E_{\bs\gamma}\) under the
parametrization~\eqref{eq:ifgf_cov}
\begin{equation}\label{eq:box_cone_seg}
    C_{\bs\gamma}(\bsx^{S}) :=
\{ \bsx = x(s,\theta,\phi;\bsx^{S}) : ((s,\theta,\phi)\in E_{\bs\gamma})\},
\end{equation}
a \emph{cone segment} (see Fig.~\ref{fig:cone_struct_2d} for an
analogous two-dimensional illustration for cone-structure for a single box.) 
We point out that by definition,
\(\mcE = \cup_{\bs\gamma\in K_{C}} C_{\bs\gamma}(\bsx^{S})\) and
\(C_{\bs\gamma}(\bsx^{S})\cap C_{\bs\gamma'}(\bsx^{S}) = \emptyset\)
for \(\bs\gamma\neq\bs\gamma'\).

In each cone segment \(C_{\bs\gamma}\), the IFGF algorithm uses 1D
Chebyshev interpolation strategy in each of the radial and angular
variables with a fixed numbers \(P_{s}\) and \(P_{a}\) of
interpolation points along the radial and the angular variables,
respectively. For a function \(f:[-1,1]\rightarrow \mb{C},\) the 1D
Chebyshev interpolating polynomial of degree \((n-1)\) is given by 
\begin{equation}
    I^{n}f(x) = \sum_{i=0}^{n-1} \widetilde{w}_{i}T_{i}(x),\quad x\in [-1,1],
\end{equation}
where \(T_{i}\) is the Chebyshev polynomial of degree \(i\) and the
weights \(\widetilde{w}_{i}\) are given by 
\begin{equation}
    \widetilde{w}_{i}=\frac{2}{n}{\sum}_{j=0}^{n-1}c_{j}f(x_{j})T_{i}(x_{j}),\quad
\mathrm{ for }\; i=0,\dots,n-1
\end{equation}
with \(c_{0}=0.5\) and \(c_{i}=1\) for \(i=1,\dots,n-1\), and
the points \(x_{j}\) are given by~\eqref{eq:cheby_nodes} for \(n=N_{C}\).
  
\subsection{IFGF Multilevel Recursive Interpolation.}\label{sec:ifgf_multi_lvl}
To achieve the desired acceleration with a \(\mc{O}(N \log N)\) computational cost, the
proposed method utilizes the multi-level IFGF recursive-interpolation
strategy introduced in~\cite{ifgf_jcp} for the Helmholtz problem. The
multi-level IFGF method implements the 
interpolation strategy described for a single box in
Section~\ref{sec:ifgf_single_box} in a recursive 
manner using larger and larger boxes, where the 
contributions of the larger boxes are in turn evaluated by
interpolation and accumulation of the contributions from the smaller
boxes.

\begin{figure}[ht]
  \centering
  \includegraphics[scale=0.42]{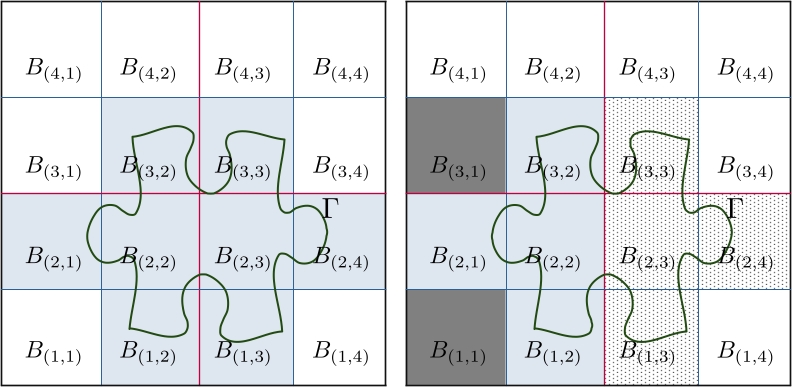}
  \caption{\footnotesize{Illustration for box-hierarchy in 2D for \(D=3\). The level
   \(1\) box (the largest box) is always relevant as it contains the set
\(\Gamma_{N}\). Here, the all four level-\(2\) boxes (divided by the red
lines) are relevant as well. On level \(3\)--On left: only the light shaded
boxes on the left image above are
relevant as the other boxes do not intersect with the surface
\(\Gamma\). On right: for the box \(B_{(2,1)}\), only the light shaded boxes
are in \(\mc{N}B_{(2,1)}\) whereas the dot-pattern boxes are in the
cousin box set \(\mc{M}B_{(2,1)}\). Note that the dark shaded boxes
\(B_{(1,1)}\) and \(B_{(3,1)}\) not being relevant are not considered
neighbors of \(B_{(2,1)}\)~\eqref{eq:box_neighbor}.}\label{fig:box_ills_2d}} 
\end{figure}

The multi-level IFGF method starts by selecting a single cubic box
\(B(\bsx_{(1,1,1)},H_{1}) \supset \Gamma_{N} \)
(see~Fig.~\ref{fig:box_ills_2d}) containing all 
the source points. Then, starting from this topmost
level-\(d=1\) box \(B(\bsx_{(1,1,1)},H_{1})\), the IFGF method
creates the level-\(d\) boxes for \(d=2,\dots,D\), by partitioning
each level-\((d-1)\) boxes into eight disjoint, equisized
\emph{``children''} boxes 
\(B^{d}_{\bs{\mr k}}=B(\bsx_{\bs{k}},H_{d})\) of side \(H_{d}=H_{d-1}/2\), and centered at
\[
  \bsx^d_{\bs{\mr k}} := \bsx^1_{\bs{1}} - \frac{H_{1}}{2^{}} \bs{1} +
\frac{H_d}{2^{}} (2\bs{\mr{k}}-\bs{1}),
 \; 
(\bs{\mr k}\in \mr K^d_{B}=\{1,\dots,2^{d-1}\}^3).
\]
Each box \(B^{d}_{\bs{\mr k}}\) for \(d=2,\dots,D\) is contained in a
 ``parent box'' on level-(\(d-1\)), which is denoted by
\(\mc{P}B^{d}_{\bs{\mr k}}\). The value of \(D\) is chosen so that the
side \(H_{D}\), corresponding to the side length of a box on the lowest level, satisfies the condition \(H_{D} \leq 0.25\lambda\),
where \(\lambda\) is the wavelength. It is only necessary to
keep track of the \emph{``relevant boxes''}, that is, the boxes
\(B^{d}_{\bs{\mr k}}\) for which \(B^{d}_{\bs{\mr k}}\cap
\Gamma_{N}\neq\emptyset\). The set of all relevant boxes of level
\(d\) (\(d=1,\dots,D\)) is denoted by 
\begin{equation}\label{eq:rele_boxes}
    \mc{R}^{d}_{B}:=\{B^{d}_{\bs{\mr k}}: B^{d}_{\bs{\mr k}}\cap
\Gamma_{N}\neq\emptyset; \bs{\mr k}\in \mr K^d_{B}\}.
\end{equation}  
The totality of the above box hierarchy is stored as a linear octree. 

Utilizing~\eqref{eq:box_contr_factor} and~\eqref{eq:box_analfac}, the
field generated by the sources within a relevant box \(B^{d}_{\bmk}\) are
given by
\begin{equation}\label{eq:bkd_contr}
  I^{d}_{\bmk}(\bsx) =
  G(\bsx,\bsx^{d}_{\bmk}) F^{d}_{\bmk}(\bsx);\;
  F^{d}_{\bmk}(\bsx) =\!\!\!\!
\displaystyle\sum_{\bsy\in B^{d}_{\bmk}\cap\Gamma_{N}}\!\!\!\!\!\!\!
 a(\bsy) g^{d}_{\bmk}(\bsx,\bsy).
\end{equation}
Moreover, using~\eqref{eq:ifgf_cov}
\begin{equation}\label{eq:bkd_contr_ana_fac}
   F^{d}_{\bmk}(\bsx) = F^{d}_{\bmk}(s,\theta,\phi) 
   = \displaystyle\sum_{\bsy\in B^{d}_{\bs{\mr{k}}}\cap\Gamma_{N}}\!\!\!\!\!\!\!
 a(\bsy) g^{d}_{\bmk}(s,\theta,\phi,\bsy)
\end{equation}
where the spherical coordinate \((s,\theta,\phi)\) is centered at
\(\bsx^{d}_{\bs{\mr k}}\). As discussed in what 
follows, the multilevel recursive interpolation strategy relies on the
application of the single box interpolation strategy to evaluate the
analytic factor \(F^{d}_{\bs{\mr k}}(s,\theta,\phi)\) for each one of
the relevant boxes starting at level-\(D\), and then
iteratively proceeding to level-\(3\); at this level for a given target
point \(\bsx\in B^{d}_{\bs{\mr{k}}}\) contributions from all the
non-neighbor boxes are accumulated. The contribution from the sources
within the \emph{``neighbor boxes''} (defined in what follows) are
computed using the CBIE method as discussed earlier in 
Sections~\ref{sec:non_sing} and~\ref{sec:rect_pol}. 

To facilitate the recursive interpolation strategy, several additional
concepts are required. Following~\cite{ifgf_jcp}, we define for a given
level-\(d\) box \(B^{d}_{\bs{\mr k}}\), the set \(\mc{N}
B^{d}_{\bs{\mr k}}\) of all level-\(d\) boxes that are neighbors of
\(B^{d}_{\bs{\mr k}}\) (that is a box that shares a side with the box
\(B^{d}_{\bs{\mr k}}\)) and, the set \(\mc{M} B^{d}_{\bs{\mr k}}\) of
all level-\(d\) boxes that are cousins of \(B^{d}_{\bs{\mr k}}\) (the
non-neighboring boxes which are children of a neighbor of the 
parent box \(\mc{P}B^{d}_{\bmk}\)). Similarly, the set \(\mc{U} B^{d}_{\bs{\mr k}}\)
of \emph{neighbor points} and the set \(\mc{V} B^{d}_{\bs{\mr k}}\) of
\emph{cousin points} of \(B^{d}_{\bs{\mr k}}\) are defined as the set
of all points in \(\Gamma_{N}\) that are contained in the neighbor and
the cousin boxes of \(B^{d}_{\bs{\mr k}}\), respectively. We thus have 
\begin{align}
  \mc{N} B^{d}_{\bs{\mr k}} &:=
  \{B^{d}_{\bs{\mr j}}\in\mc{R}^{d}_{B} : \lVert \bs{\mr j}-\bs{\mr k} \rVert \leq 1 \}\label{eq:box_neighbor} \\
  \mc{M} B^{d}_{\bs{\mr k}} &:=
  \{ B^{d}_{\bs{\mr j}}\in\mc{R}_{B}: B^{d}_{\bs{\mr j}}\notin \mc{N}B^{d}_{\bs{\mr k}} \wedge \mc{P} B^{d}_{\bs{\mr j}}\in \mc{NP} B^{d}_{\bs{\mr k}}\}\label{eq:box_cousin} \\
  \mc{U} B^{d}_{\bs{\mr k}} &:=
  \left( \bigcup_{B\in\mc{N}B^{d}_{\bs{\mr k}}} B \right) \cap \Gamma_{N}\label{eq:neighbor_points} \\
  \mc{V} B^{d}_{\bs{\mr k}} &:=
  \left( \bigcup_{B\in\mc{M}B^{d}_{\bs{\mr k}}} B \right)\cap \Gamma_{N}\label{eq:cousin_points}
\end{align}
Figure~\ref{fig:box_ills_2d} (see captions thereof) illustrates an
analogous two-dimensional box-hierarchy for \(D = 3\) and the
aforementioned concepts, namely, relevant boxes, neighbors and cousins.

As alluded to in the single box interpolation discussion in
Section~\ref{sec:ifgf_single_box}, the interpolation of the
analytic factor \(F^{d}_{\bs{\mr k}}(s,\theta,\phi)\) associated with 
a level-\(d\) relevant box \(B^{d}_{\bs{\mr k}}\) relies on
level-\(d\) cone segments, which are denoted by
\begin{equation}\label{eq:cone_hier}
    C^{d}_{\bs{\mr k};\bs\gamma} = C_{\bs\gamma}(\bsx^{d}_{\bs{\mr
k}}), \; \bs\gamma\in K^{d}_{C}
\end{equation}
where
\(K^{d}_{C}:= \{1,\dots,n_{s;d}\} \times \{1,\dots,n_{a;d}\}\times
\{1,\dots,2 n_{a;d}\}\) 
for certain selection of integers \(n_{s}=n_{s;d}\) and
\(n_{a}=n_{a;d}\) appropriately chosen
to ensure the desired accuracy (see
Remark~\ref{rmk:cone_struct_size}). Analogous to the relevant 
boxes, the IFGF method also uses the concept of \emph{``relevant cone
segments''}. Denoting the set of interpolation points within the cone
segments \(C^{d}_{\bs{\mr k};\bs\gamma}\) by \(\mc{X} C^{d}_{\bs{\mr
k};\bs\gamma}\); the set \(\mc{R}_{C}B^{d}_{\bs{\mr k}}\) of the
cone segments relevant to the box \(B^{d}_{\bs{\mr k}}\) are defined
recursively starting from level \(d=3\) as follows: 
\begin{align*}\label{eq:rele_cone}
  \mc{R}_{C}B^{d}_{\bmk} &:= \emptyset \; \mathrm{ for }\; d=1,2,\,\mathrm{and,}\\
  \mc{R}_{C}B^{d}_{\bmk} &:= \Big\{ C^{d}_{\bmk;\bg}:
  C^{d}_{\bmk;\bg}\cap\mc{V}B^{d}_{\bmk}\neq\emptyset\,
  \mathrm{ or }\, C^{d}_{\bmk;\bg}\cap
\Big(\hspace*{-5mm}
  \bigcup_{C\in\mc{R}\mc{P}B^{d}_{\bmk}}\hspace*{-4mm} \mc{X}C \Big)\neq\emptyset
  \Big\}
\end{align*}
for \(\bmk\in \mr K^d_{B}\) and \(d\geq 3\). Thus, a level-\(d\) cone
segment \(C^d_{\bmk;\bs\gamma}\) is recursively (starting from \(d=3\)
to \(d=D\)) defined to be relevant to a box \(B^d_{\bmk}\) if either, 
(1) it includes a surface discretization point in a cousin of
\(B^d_{\bmk}\),  or if, (2) it includes an interpolation point of a relevant cone
segment associated with the parent box \(\mathcal{P}B^d_{\bmk}\) of
\(B^d_{\bmk}\). The set of all relevant level-\(d\) cone-segments is
given by 
\begin{equation}\label{eq:rele_cone_d}
 \mc{R}_{C}^{d} \coloneqq \{C^{d}_{\bmk;\bs{\gamma}}\in 
\mc{R}_{C}B^{d}_{\bmk}: \bs{\gamma}\in K^{d}_{C}, \bmk\in K^{d}_{B}
~\mathrm{and}~ B^{d}_{\bmk}\in \mc{R}^{d}_{B} \}.
\end{equation}

\begin{remark}[\textbf{Choosing cone-structure sizes.}]\label{rmk:cone_struct_size}
\normalfont{In view of Theorems 1 and 2 in~\cite{ifgf_jcp}, which in
particular imply that the use of fixed numbers 
\(n_{s}\) and \(n_{a}\) give rise to essentially constant cone-segment
interpolation errors for all \(H\)-side boxes satisfying \(\kappa
H<1\), and the pertaining discussion in \cite[Section 3.3.1]{ifgf_jcp},
the IFGF algorithm chooses the values of \(n_{s;d}\) and \(n_{a;d}\) as
described in what follows. Starting from certain selection of
\(n_{s;(D+1)}\) and  \(n_{a;(D+1)}\) as the initial values, the
method selects the value of \(n_{s;d}\) and  \(n_{a;d}\) for
\(d=D,D-1,\dots,3, \) according to the rule
\(n_{s;d}=a_{d}n_{s;(d-1)}\) and  \(n_{a;d}=a_{d}n_{a;(d-1)}\)
where \(a_{d}=1\) if \(\kappa H_{d} \leq 1/2\), and \(a_{d}=2\) in the
complementary case. In keeping with \cite{ifgf_jcp}, in each
cone-segment we have used the values \(P_{s}=3\) and \(P_{a}=5\) for
all the numerical experiments presented in this work.}
\end{remark}

As mentioned above, the IFGF method evaluates the discrete integrator
by commingling the effect of large numbers of sources into a small
number of interpolation parameters; which is possible, in particular,
owing to the analytic property and slow-oscillatory character of the
analytic factor \(g^{d}_{\bmk}\).
 We note that at the level \(D\), the contribution
of the sources within a box \(B^{D}_{\bmk}\) at the relevant
cone-segment interpolation points are produced by directly evaluating
the sum \(F^{D}_{\bmk}\). A summary of the algorithm is
presented below.

\begin{enumerate}
\item Determine the sets of relevant boxes
\(\mc{R}_{B}\)~\eqref{eq:rele_boxes} and of relevant cone-segments
\(\mc{R}^{d}_{C}\)~\eqref{eq:rele_cone_d} for all \(d=1, \dots, D\).
\item Direct evaluations on level \(D\).\\
   --For every D-level box evaluate the field \(I^{D}_{\bmk}\) at all
the neighboring target points \(\bsx\in \mc{U} B^{D}_{\bmk}\) using the CBIE method.\\
   --For every D-level relevant box \(B^{D}_{\bmk}\in \mc{R}_{B}\)
evaluate the analytic factor \(F^{D}_{\bmk}\) at all the interpolation
points \(\bsx\in\mc{X} C^{D}_{\bs{\bmk};\bs\gamma}\) of the relevant cone
segments \(C^{D}_{\bs{\bmk};\bs\gamma}\in \mc{R}_{C}B^{D}_{\bmk}\) by
direct evaluation of the sum~\eqref{eq:bkd_contr_ana_fac}.
\item Interpolation for \(d = D,\dots,3\).\\
   --For every relevant box \(B^d_{\bmk}\) evaluate the field
\(I^d_{\bmk}(\bsx)\) at every surface discretization point \(\bsx\in \mc{V} B^{d}_{\bs{\mr k}}\)
by interpolation of the
analytic factor \(F^d_{\bmk}\) and multiplication by the centered
factor \(G(\bsx,\bsx^d_{\bmk})\).\\    
   --For every level-\(d\) relevant box \(B^d_{\bmk}\) determine the parent box
\(B^{d-1}_{\bs{\mr j}} = \mc{P} B^d_{\bmk}\) and, obtain the analytic
factor \(F^{d-1}_{\bs{\mr j}}\) at all level-\((d-1)\) interpolation
points corresponding to \(B^{d-1}_{\bs{\mr j}}\) by first,
interpolation of \(F^d_{\bmk}\) and then re-centering by multiplying
the smooth factor \(G(\bsx,\bsx^d_{\bmk})/ G(\bsx,\bsx^{d-1}_{\bs{\mr
j}})\), and finally, accumulating the contributions from all the children
boxes of \(B^{d-1}_{\bs{\mr j}}\).
\end{enumerate}
For a more detailed description of the IFGF acceleration strategy,
including an analysis of its computational complexity, we
refer to~\cite{ifgf_jcp}. 


\subsection{Local Corrections to IFGF Sum.}\label{sec:local_corr}
\begin{figure}[ht]
  \centering
  \includegraphics[scale=0.5]{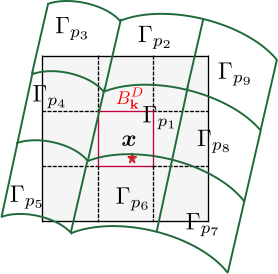}
  \caption{\footnotesize{Illustration for local correction in 2D.}\label{fig:local_corr}}
\end{figure}

For evaluation of \(I^{D}_{\bmk}\) at the points \(\bsx_{\ell}\in \mc{U}B^{D}_{\bmk}\) the 
method relies on the CBIE part of the method. As
discussed in Section~\ref{sec:surf_discr}, the method uses a
non-overlapping patch decomposition of the surface 
\(\Gamma\), and certain changes of variables on a refined grid to
resolve the singular behavior of the kernel when the source and the
target points are in close proximity. It becomes necessary to
compute the total contribution from all the source points within a
patch as the quadrature is target point specific. Clearly, the surface
decomposition does not conform to the 3D 
boxes, and some of the patches covering the neighbor set
\(\mc{N}B^{D}_{\bmk}\) may intersect the cousin set \(\mc{M}B^{D}_{\bmk}\). 
Note that if any of these intersecting patches, say, \(\Gamma_{p}\)
is singular for a point \(\bsx\in B^{D}_{\bmk}\), that is,
\(p\in\mathscr{S}_{\ell}\) (see~\eqref{eq:patch_sing}) then some
contribution from these patches, for example, the region  
\(\Gamma_{p_{6}}\setminus \mc{N}B^{D}_{\bmk}\) as per the illustration
 in Figure~\ref{fig:local_corr}, is counted both by the CBIE
 and the IFGF part of the algorithm, which requires a
correction. We note that these corrections can be applied at the
beginning of a run by subtracting the product of Fej\'{e}r quadrature
weight and the kernel value from the precomputed singular weights for
all such instances.

\subsection{IFGF for EM Problem.}\label{sec:ifgf_vec}
In Section~\ref{sec:ifgf_multi_lvl}, we have discussed the recursive
IFGF interpolation strategy for a single density \(\phi\). For the EM problem,
as mentioned at the beginning of Section~\ref{sec:numer_meth}, we must
carry out the interpolation strategy for eight scalar quantities, and we
also need to compute the normal derivative of the single layer (SL) potentials as well.\\
\textbf{Computing Normal Derivative.} To compute the normal derivative
of the SL potentials, we use finite differences, in
particular, we use the (first order) forward difference formula 
\begin{equation}\label{eq:der_fd}
D[\phi](\bsx) = \frac{S[\phi](\widetilde{\bsx}_{\bsn}) - S[\phi](\bsx)}{|\bsn\cdot\bs{h}|}
\end{equation}
where \(\widetilde{\bsx}_{\bsn}=\bsx + \bsn\cdot\bs{h}\) and
\(\bs{h} = (h,h,h)^{t}\) with one-dimensional step-size 
\(h=|\bsn\cdot\bs{h}|=10^{-6}\). 
Hence, to compute the normal derivative we need to evaluate the
interpolation at one extra point for each density for which the normal
derivative is required. To obtain the value
\(S[\phi](\widetilde{\bsx_{n}})\) for a surface discretization point
\(\bsx\in B^{d}_{\bmk}\), the method uses the same interpolation
procedure. In particular, to obtain the contribution to the normal
derivative \(D[\phi](\bsx)\) from a cousin box 
\(B\in\mc{M}B^{d}_{\bmk}\), the method uses the relevant cone-segment
containing the point \(\bsx\) to interpolate at the point
\(\widetilde{\bsx}_{\bsn}\) in addition to the point \(\bsx\)
itself. Then the formula~\eqref{eq:der_fd} is used to approximate the
contribution to the normal derivative at \(\bsx\) coming from a
cousin box. 
Note that for a given point \(\bsx\in\Gamma_{N}\), it is possible for the
point \(\widetilde{\bsx}_{\bsn}\) to lie in a different
cone-segment, which leads to extrapolation instead of interpolation. The
point \(\widetilde{\bsx}_{\bsn}\) may also lie within the neighbor
boxes. Both situations, due to the small finite-difference
step-size $h$, do not affect the accuracy as demonstrated by the
numerical experiments presented in Section~\ref{sec:numer_exp}.

\textbf{IFGF for Multiple Densities: Speed-vs-Memory Trade-off.}\\
As mentioned in the beginning of this section, in order to efficiently
approximate the integral operators associated with~\eqref{eq:n_muller},
we need to use the IFGF interpolation strategy for eight densities and
two kernel values arising from two different domains. The
straightforward approach would be to use the IFGF algorithm once for
each combination of a density and a kernel, which here onward we call
\(\mathrm{IFGF}_{Seq}\). In this approach, we are
required to traverse the box-octree
a total number of sixteen times. A different
approach would be to interpolate the integrals for 
all the densities and the kernels at one go, which we refer to as
\(\mathrm{IFGF}_{Vec}\). In the later approach, we need to traverse
the box-octree only once. More importantly, in the second, the kernel values for a given pair of source and
target points only need to be computed once for all the densities, leading to a further
speedup of the overall IFGF strategy in the context of the EM problem compared to the
\(\mathrm{IFGF}_{Seq}\) approach (see
Fig.~\ref{fig:time_ifgf_seq}). We point out that in the \(\mathrm{IFGF}_{Seq}\)  
approach, we need to store the interpolation weights only for one 
combination of a density and a kernel at any time whereas in the \(\mathrm{IFGF}_{Vec}\)
approach, we need to store the interpolation weights for all sixteen
possible combinations of the densities and both the kernels simultaneously, which
requires sixteen times the memory for storing the interpolation
weights in comparison to the first approach.

For this work, we chose the second option, which leverages uses the same
IFGF algorithm, but for an array of densities and an array of
kernels, to save on the run time (see Remark~\ref{rmk:run_memory}). In
particular, the input weights to the IFGF algorithm for the EM problem
can be written as \(W=(a_{k,\ell})_{8\times N}\). Denoting by \(I_{k,n}\)
the sum for \(k\)-th density and \(n\)-th kernel, and denoting the two
kernels by \(G_{1}\) and \(G_{2}\), the EM version of 
equation~\eqref{eq:ifgf_sum} can be written as  
\begin{equation}\label{eq:ifgf_sum_em}
\begin{split}
  \begin{bmatrix}
    I_{1,1}(x_\ell) & I_{1,2}(x_\ell)\\
    \vdots & \vdots\\
    I_{8,1}(x_\ell) & I_{8,2}(x_\ell)
  \end{bmatrix}&
  \\
    := \sum_{\substack{m=1 \\ m\neq\ell}}^{N}
     \begin{bmatrix}
       a_{1,m}\\
       \vdots\\
       a_{8,m}
     \end{bmatrix}
  &
   \begin{bmatrix}
      G_{1}(\bsx_\ell,\bsx_m) & G_{2}(\bsx_\ell,\bsx_m)
    \end{bmatrix},
  \end{split}
\end{equation}
for \(\ell=1,\dots,N\). The evaluation of \eqref{eq:ifgf_sum_em} then
follows the steps of the IFGF algorithm prescribed in
Section~\ref{sec:ifgf_multi_lvl}. This concludes the discussion on the
numerical scheme of the proposed methodology.

\begin{remark}[\textbf{On selection of IFGF approach}]\label{rmk:run_memory}
\normalfont{One can choose to approximate the discrete sum for a
specific number of densities or kernels at a time. For instance, we can evaluate the sum
for all eight densities, but for a single kernel at a time; this reduces the 
memory requirement of the interpolation weights exactly by half, and
(in our experiments) shows \(\approx 1.5\times\) slower run time
compared to \(\mathrm{IFGF}_{Vec}\). Moreover, the memory required by
the interpolation weights can easily be computed once we know the
number of relevant cone segments; and one can choose an
approach at \emph{run time} based on the available memory.}
\end{remark}

\section{Numerical Experiments}\label{sec:numer_exp}
In this section, we present results from several numerical
experiments. In all the examples presented, equal
numbers of mesh points are used in \(u\) and \(v\) variables to
discretize each patch. For the IFGF interpolation,
\(P_{s}=3\), and \(P_{a}=5\) values were used in the radial and
angular directions, respectively. In addition, we have used the values
\(n_{s;D+1}=1\) and \(n_{a;D+1}=2\) (see~\eqref{eq:cone_hier}) to
structure the cone-hierarchy. 
For the incident field, we considered the
planewave \(\bsE^\mr{inc} = \exp(-i\kappa_{e}z)\hat\bsx\) in all cases, except in
Example~\ref{sec:exm_splitter}, where an electric dipole was used
instead. We recall that \(N=P N^{2}_{C}\) denotes the discretization size~\eqref{eq:mesh_pts_l} (the total
number of unknowns is \(4N\)), \(P\) denotes the number of
decomposing patches~\eqref{eq:surf_patch}, and \(N_{C}\) denotes the number of points per
patch in one variable~\eqref{eq:cheby_nodes}.
We start the experiments with a forward map computation
for a spherical geometry to study the time and accuracy of the
CBIE-IFGF method. In addition, for the spherical geometry, 
we study the scattering simulation, as in this case the Mie
series solution is available to compute the accuracy of the proposed
method. 
Next, in order to demonstrate the applicability and the performance of
the CBIE-IFGF method for arbitrarily 
shaped geometries, we also consider several different CAD models, namely, a glider model in
Example~\ref{sec:exm_glider}, a hummingbird model in Example~\ref{sec:exm_hummingbird}, and 
a \textit{Gmsh}-rendered nanophotonic power splitter in Example~\ref{sec:exm_splitter}. Run times of
one forward map computation of the accelerated CBIE-IFGF method are
provided for all examples and compared against that of the
unaccelerated CBIE method.
 Note
that the unaccelerated CBIE method implementation includes certain additional optimizations (e.g., reduced number of
kernel evaluations due to shared kernels between densities) that are unavailable for the CBIE-IFGF method.
The simulations for all the numerical results presented in this work
were run on \(128\) cores 
of an Ubuntu server equipped with two AMD EPYC \(7763\) \(64\)-Core
Processors (CPU speed: max \(3.5\) GHz and min \(1.5\) GHz) and a
total of \(1\)TB available RAM. The GMRES residual tolerance was set to \(10^{-4}\).

\subsection{Forward Map Computation\label{sec:exm_fmap}}
As a first experiment, we study the timings and the accuracy
of the proposed (accelerated) CBIE-IFGF forward map (FM)---namely, the
action of the discretized version of the operators in the l.h.s.~of
\eqref{eq:n_muller} on a given set of densities. For this experiment, we consider
\(\Omega_{i}\) to be the 
unit sphere centered at the origin; \(\Omega_{e}\) being the exterior
\(\mathbb{R}^{3}\setminus\Omega_{i}\) with \(\epsilon_{e}= 1,
\epsilon_{i}=2.25\), and
\(\mu_{e}=\mu_{i}=1\). The refractive index of the interior medium is \(1.5\).
For this simulation, we have used the value
\(N_{C}=16\) (the smallest discretization in Table~\ref{tab:exm_fm}
has \(P=24\) patches, and \(P\) increases by a factor of four in each
subsequent finer discretizations; this choice provides an accuracy of
order \(10^{-5}\) or better in
the solution obtained using the CBIE method) to get an estimate of the accuracy
provided by the CBIE-IFGF method. The Err.~(FM) column in
Table~\ref{tab:exm_fm} shows the \(\ell_{2}\)-difference
between the FM values obtained by the CBIE and the CBIE-IFGF methods.
It shows that the FM operator values obtained using the CBIE-IFGF
method match with that of the CBIE
method with a fixed \(10^{-4}\) accuracy, while maintaining a fixed
number of points per wavelength. (Note that the choice of the
\(N_{C}\)-value does not affect the accuracy provided by the IFGF
interpolation strategy.) 

\begin{figure}[ht!]
  \centering
\includegraphics[scale=0.2, trim={80 0 40 40},clip]{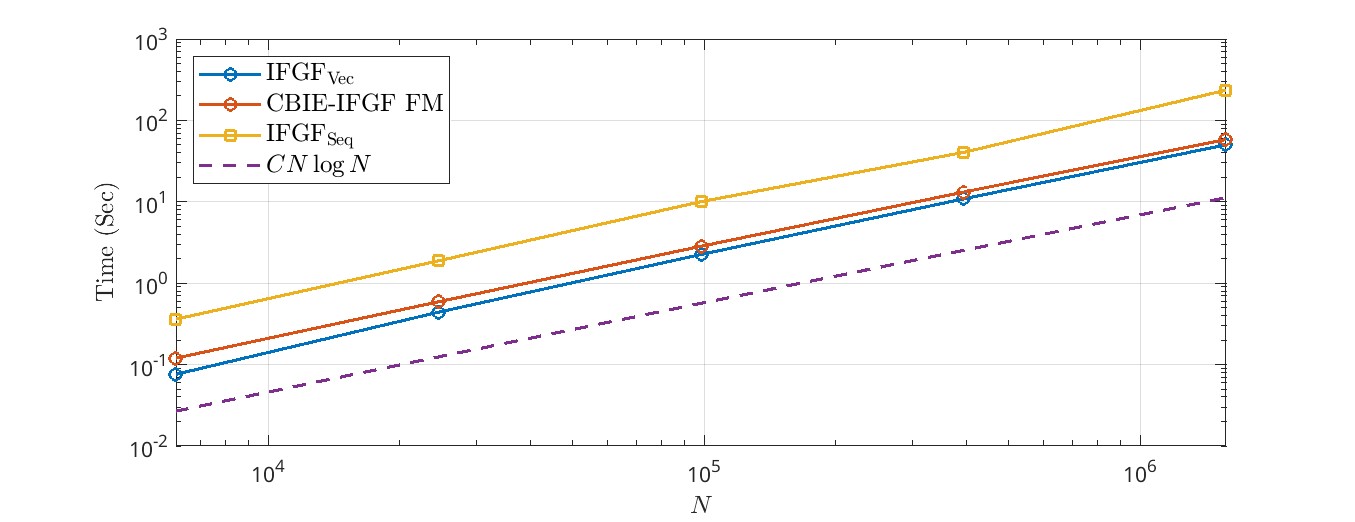}
\caption{\footnotesize{Various run times in \(\log-\log\) scale of are
presented. The legends ``\(\mathrm{IFGF}_{Seq}\)'' and
``\(\mathrm{IFGF}_{Vec}\)'' present the run time of the corresponding IFGF
approaches. The ``CBIE-IFGF FM'' legend displays the
required run time for one FM computation of the CBIE-IFGF 
method.}\label{fig:time_ifgf_seq}} 
\end{figure}

The \(T_{C}\) and the \(T_{CI}\) columns in Table~\ref{tab:exm_fm}
represent the runtime of the FM 
of the CBIE and the CBIE-IFGF methods, respectively; whereas the column 
\(T_{IFGF}\) shows the required runtime of only the \(\mathrm{IFGF}_{Vec}\)
interpolation part of the algorithm. Additionally, for comparison, the column \(T_{Seq}\) lists
the runtime of only the \(\mathrm{IFGF}_{Seq}\) interpolation strategy.
The run times, as illustrated in
Fig.~\ref{fig:time_ifgf_seq}, are consistent with the \(\mc{O}(N\log
N)\) nature of the algorithm. 
We note that there is an increase in the runtime by a factor \(\approx
4\) in the \(\mathrm{IFGF}_{Seq}\) approach in comparison to the
\(\mathrm{IFGF}_{Vec}\) strategy. 

\subsection{Scattering from Spherical Geometry.}\label{sec:exm_sphere}
For the second numerical experiment, we consider scattering by the
unit sphere. The materialistic properties and the parameter
choices being same as in the previous example in~\ref{sec:exm_fmap}
except the value of \(N_{C}\). Here we consider the value \(N_{C}=12\) 
(which, as demonstrated in Table~\ref{tab:gmres_error}, is sufficient to achieve an
accuracy of the order of \(10^{-4}\)).
The Err.(CBIE) and Err.(C-I) columns in
Table~\ref{tab:gmres_error} present the
\(\ell_{2}\)-error in the computed 
field values produced by the unaccelerated CBIE method and the accelerated
CBIE-IFGF method, respectively. The error is computed against the Mie
series solution at \(1000\) points on the 
surface of the origin centered sphere with a radius of \(0.7\). In all
the cases, a \(10^{-4}\) accuracy is maintained by the CBIE-IFGF
method. 
\begin{table}[ht]
  \centering\footnotesize
  \begin{tabular}{c c c c c c c}
    \toprule
    \(N\) & \(\kappa_{i}\) & \(T_{C}\) & \(T_{CI}\) & \(T_{IFGF}\) & Err.~(FM)& \(T_{Seq}\)\\
    \midrule
    \(6144\) & \(3\pi\) & \(0.064\) & \(0.12\) & \(0.076\) & \(1.0\cdot 10^{-4}\)&\(0.361\)\\
    \(24576\) & \(6\pi\) & \(0.909\) & \(0.591\) & \(0.440\) &\(1.2\cdot 10^{-4}\)&\(1.88\)\\
    \(98304\) & \(12\pi\) & \(13.05\) & \(2.84\) & \(2.25\) &\(1.3\cdot 10^{-4}\)&\(10.03\)\\
    \(393216\) & \(24\pi\) & \(151.05\) & \(13.15\) & \(10.85\) &\(1.3\cdot 10^{-4}\)&\(40.38\)\\
    \(1572864\) & \(48\pi\) & \(2463.5\) & \(58.14\) & \(50.14\) &\(1.4\cdot 10^{-4}\)&\(234.92\)\\    
    \bottomrule
  \end{tabular}
\caption{\footnotesize{Accuracy and run time study for the
 rectangular-polar method and
 IFGF-acceleration. \(N_{C}=16\) points per patch (in one direction). \(T_{C}\)--runtime of CBIE FM,
 \(T_{CI}\)--runtime of CBIE-IFGF FM,
 \(T_{IFGF}\)--runtime of \(\mathrm{IFGF}_{Vec}\), and \(T_{Seq}\)--run time
for \(\mathrm{IFGF}_{Seq}\).}\label{tab:exm_fm}}
\end{table}

\begin{table}[ht]
  \centering\footnotesize
  \begin{tabular}{c c c c c c}
    \toprule
    \(N\) & \(\kappa_{i}\) & Iter.(CBIE) & Err.(CBIE) & Iter.(C-I) & Err.(C-I) \\
    \midrule
    \(3456\) & \(3\pi\) & \(26\) & \(8.8\cdot 10^{-5}\) & \(27\) & \(9.71\cdot 10^{-5}\)\\
    \(13824\) & \(6\pi\) & \(45\) & \(7.80\cdot 10^{-5}\) & \(48\)  & \(9.09\cdot 10^{-5}\)\\
    \(55296\) & \(12\pi\) & \(91\) & \(4.45\cdot 10^{-5}\) & \(102\) & \(1.16\cdot 10^{-4}\)\\
    \(221184\) & \(24\pi\) & \(158\) & \(1.28\cdot 10^{-4}\) & \(167\) & \(1.66\cdot 10^{-4}\)\\
    \(884736\) & \(48\pi\) & \(-\) & \(-\) & \(443\)  & \(3.82\cdot 10^{-4}\)\\    
    \bottomrule
  \end{tabular}
  \caption{\footnotesize Accuracy in the GMRES solution for the
unaccelerated CBIE method and the CBIE-IFGF method.\label{tab:gmres_error}}
\end{table}

\begin{figure}[ht]
  \centering
\includegraphics[scale=0.18,trim={0 0 0 80},clip]{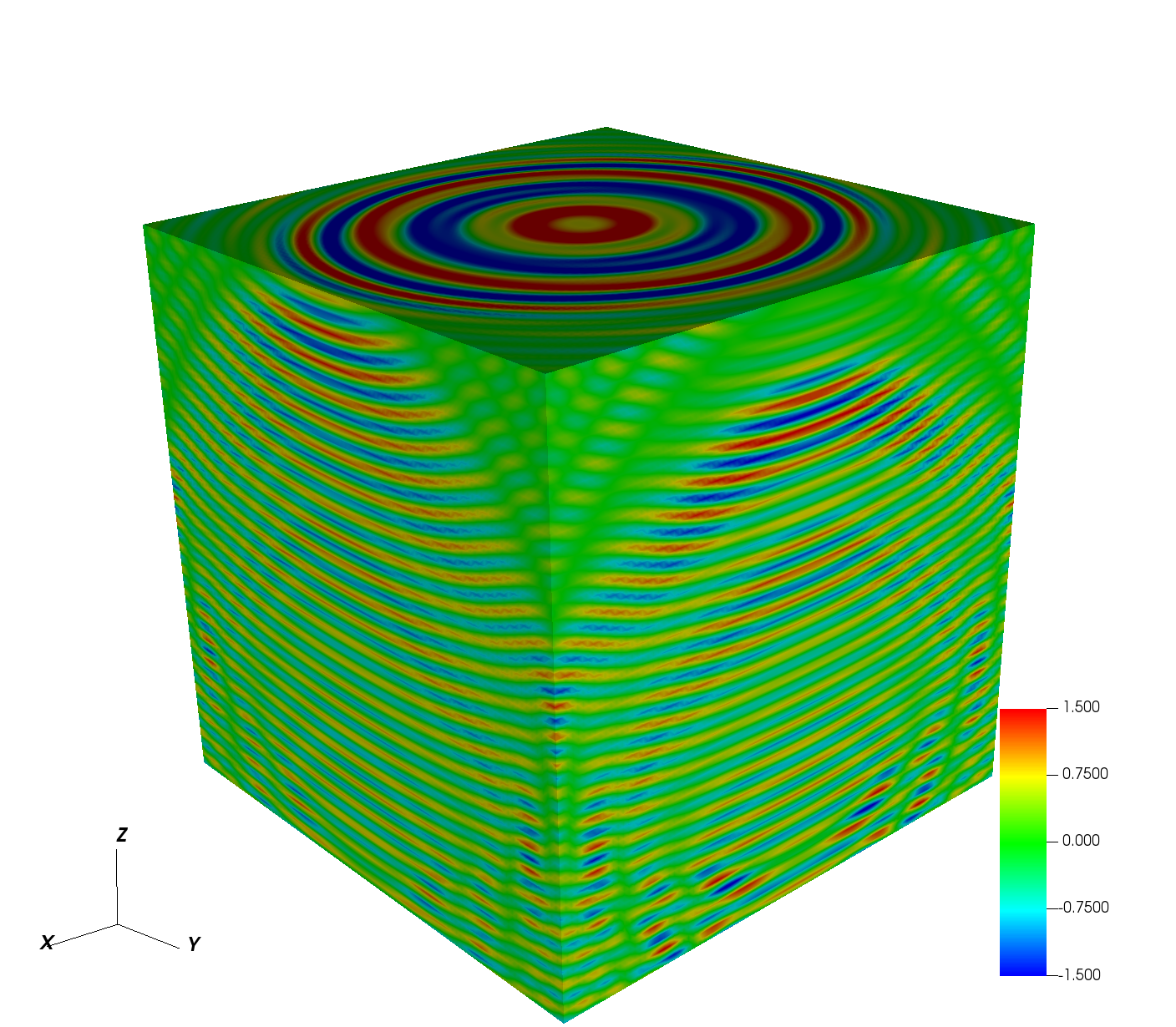}
\caption{\footnotesize Real part of the field component \(E_{x}\) for
the computed total field scattered by the unit sphere on the faces of
the origin centered cube with sides of length \(24\lambda\).
Field values are scaled to \([-1.5,1.5]\).\label{fig:sphere_field_cube}}
\end{figure}

\begin{figure}[ht]
  \centering
\includegraphics[scale=0.2,trim={0 120 0 200},clip]{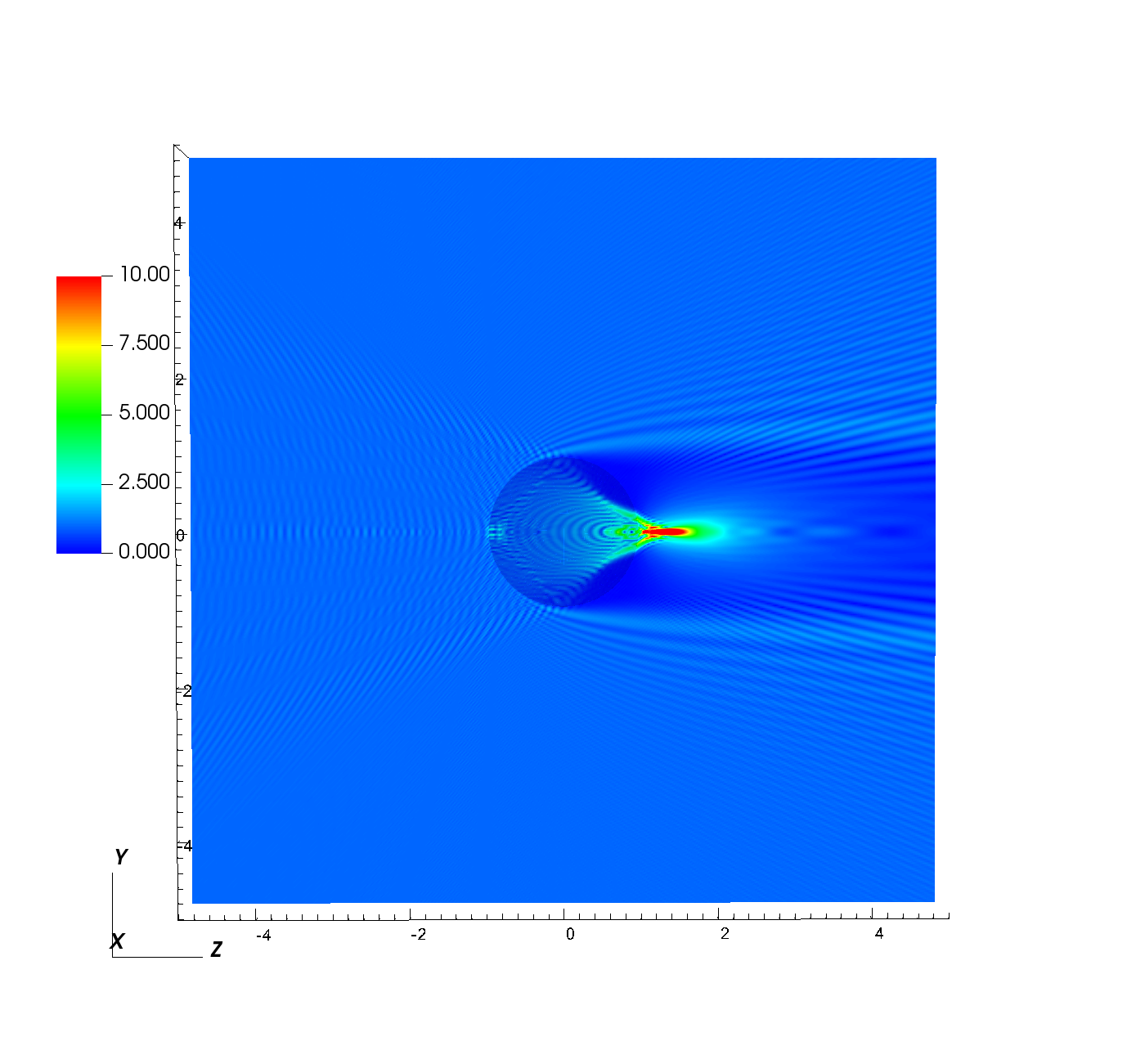}
\caption{\footnotesize Absolute value of the \(x\)-component of the
computed total field scattered by the unit sphere within the square
\([-5,5]\times [-5,5]\) on the \(yz\)-plane. A (transparent) cross-section of the
scattering sphere is inlaid in the figure. Field values are scaled to \([0,10]\).\label{fig:sphere_abs_total_yz}}
\end{figure}

For graphical presentations, we consider the values
\(\kappa_e=32\pi\) and \(\kappa_i=48\pi\).
Fig.~\ref{fig:sphere_field_cube} shows the real-part of \(E_x\)
of the computed total field scattered by the unit sphere
on the faces of the origin centered cube with sides of length
\(24\lambda\). Fig.~\ref{fig:sphere_abs_total_yz} presents the
absolute value of the \(x\)-component of the computed total field on
the \([-5,5]\times [-5,5]\) square on the \(yz\)-plane.

\subsection{Scattering From Glider CAD Model}\label{sec:exm_glider}
\begin{figure}[ht]
  \centering
\includegraphics[scale=0.25,trim={145 250 100 220},clip]{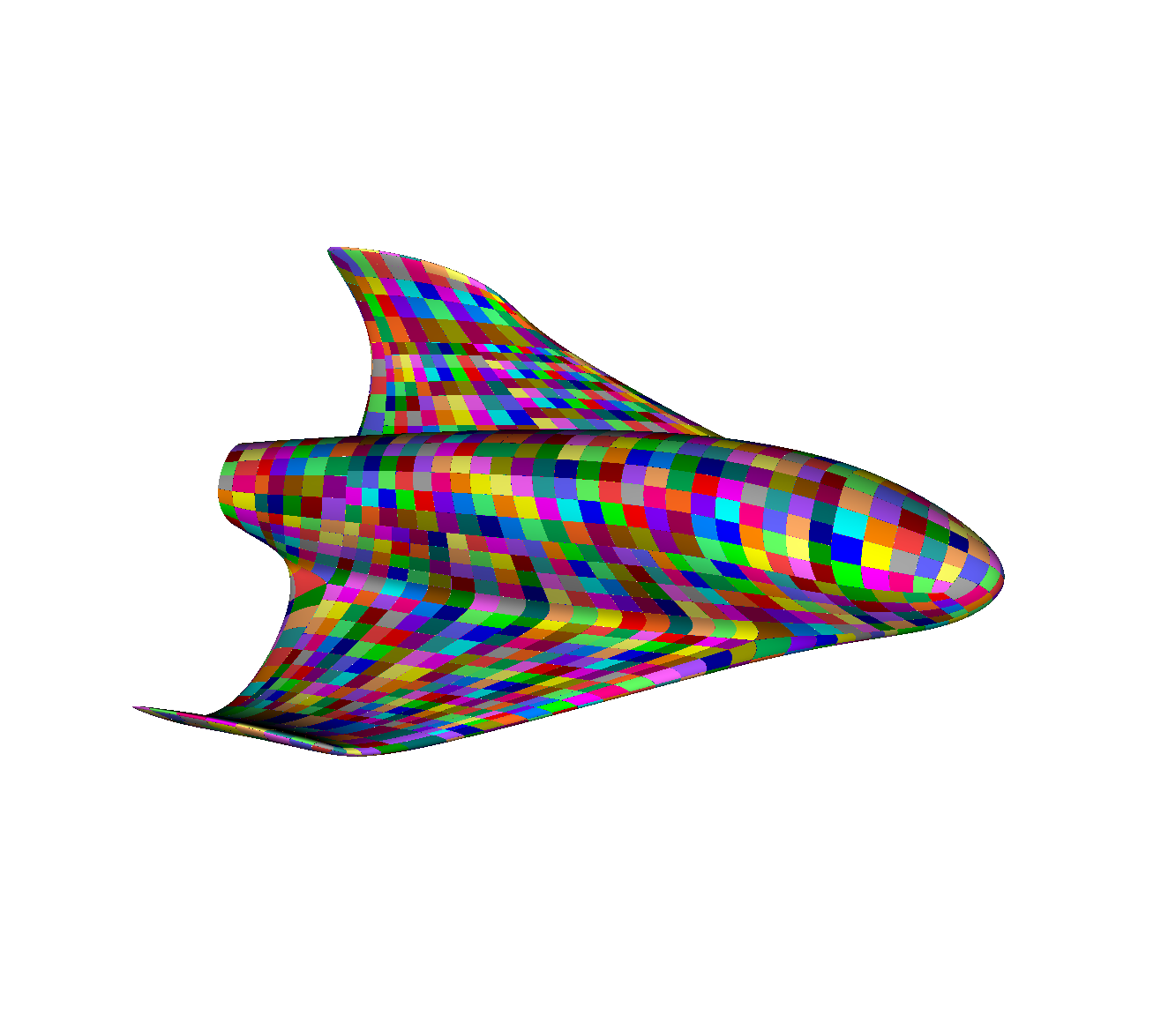}
\caption{\footnotesize{Patch configuration for the glider CAD model in
the setting of 
Example~\ref{sec:exm_glider}.}\label{fig:glider_patching}}
\end{figure}

\begin{figure}[ht]
  \centering
\begin{subfigure}[b]{1.0\linewidth}
\includegraphics[scale=0.2,trim={20 110 40 200},clip]{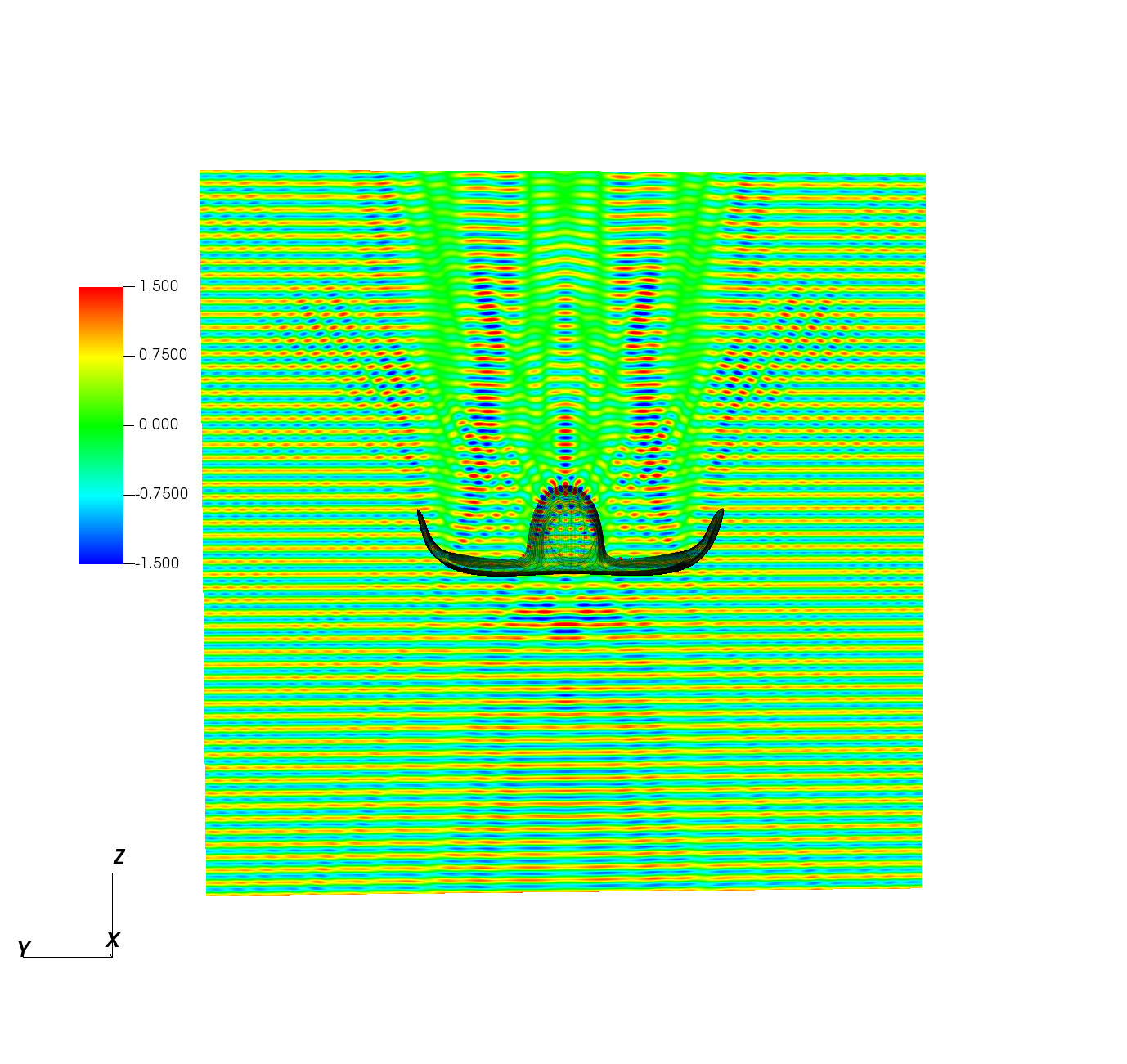}
\caption{\footnotesize{Real part of \(x\)-component of the total field
on the rectangle \([-7,7]\times [-7,7]\) on the \(yz\)-plane in the
setting of Example~\ref{sec:exm_glider}.}\label{fig:glider_ex_yz_re_rear}}
\end{subfigure}
\begin{subfigure}[b]{1.0\linewidth}
\includegraphics[scale=0.2,trim={20 110 40 100},clip]{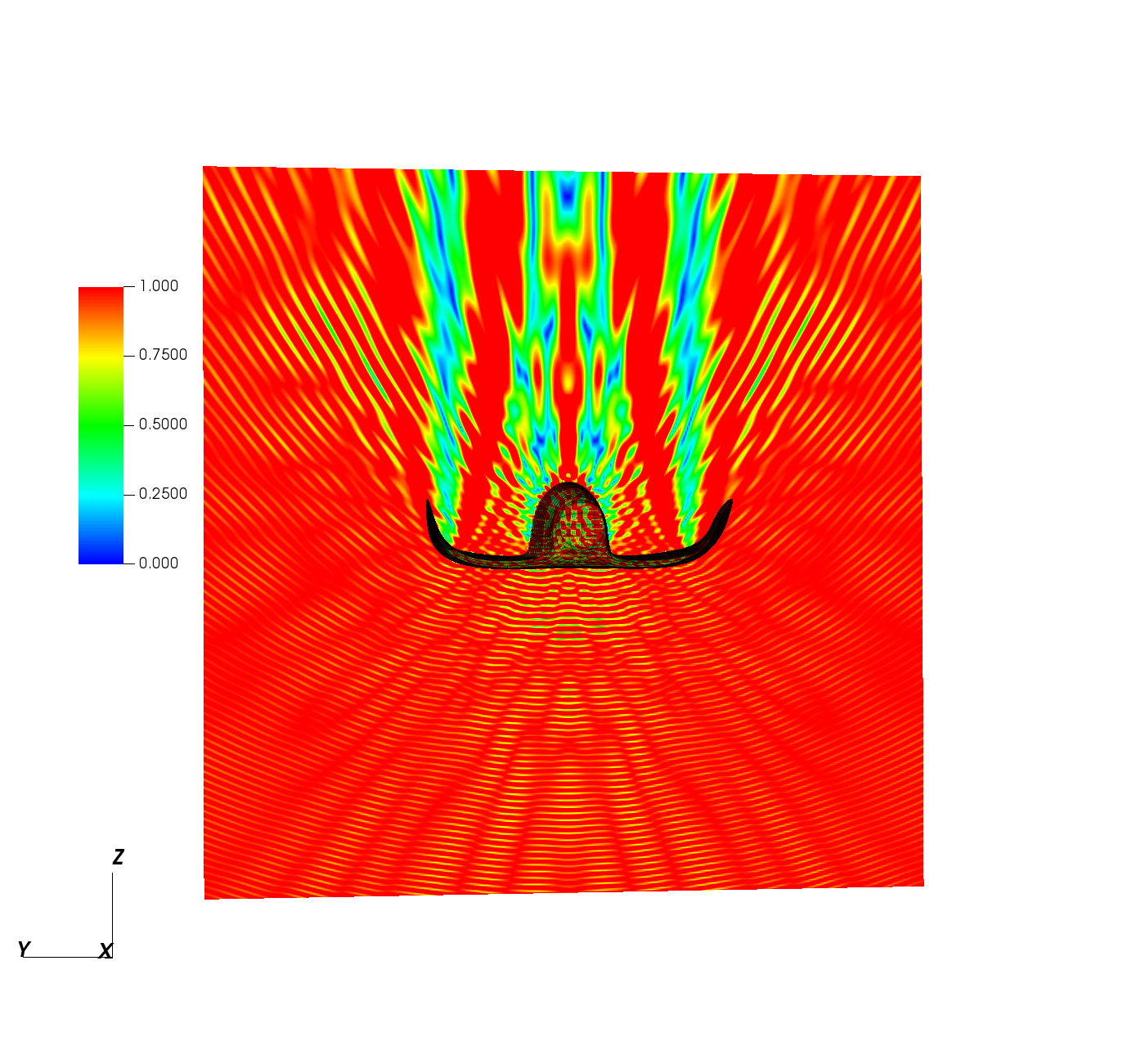}
\caption{\footnotesize{Absolute value of \(x\)-component of the total field
on the rectangle \([-7,7]\times [-7,7]\) on the \(yz\)-plane in the
setting of Example~\ref{sec:exm_glider}.}\label{fig:glider_ex_yz_abs_rear}}
\end{subfigure}
\caption{Graphical presentation for the simulation in the setting
of Example~\ref{sec:exm_glider}.\label{fig:glider_ex_yz}}
\end{figure}

For the second scattering example, we consider scattering simulation from the
glider CAD model \cite{cad_glider_src} with the following material properties: \(\epsilon_{e}=\mu_{e}=1,\kappa_{e}=25.13\), and 
\(\epsilon_{i}=2.16,\mu_{i}=1,\kappa_{i}=36.94\). The size of
the scatterer is \(45\) wavelengths in the longest dimension. 
The refractive index of the medium in \(\Omega_i\) is \(1.47\). The
discretization contains \(N\approx 320k\) points with \(P=2225\)
curvilinear patches and \(N_{C}=12\) points per patch in each variable. The
surface patch decomposition for this simulation is shown in
Fig.~\ref{fig:glider_patching}. The accelerated CBIE-IFGF solver takes
\(\approx 9\) seconds on average for one FM calculation against
\(110\) seconds by the unaccelerated CBIE
solver. Fig.~\ref{fig:glider_ex_yz} presents the real part and the absolute value of the \(x\)-component of
the computed total field.

\subsection{Scattering From Hummingbird NURBS CAD Model}\label{sec:exm_hummingbird}
\begin{figure}[ht!]
  \centering\footnotesize
\includegraphics[scale=0.2, trim={340 140 20 100},clip]{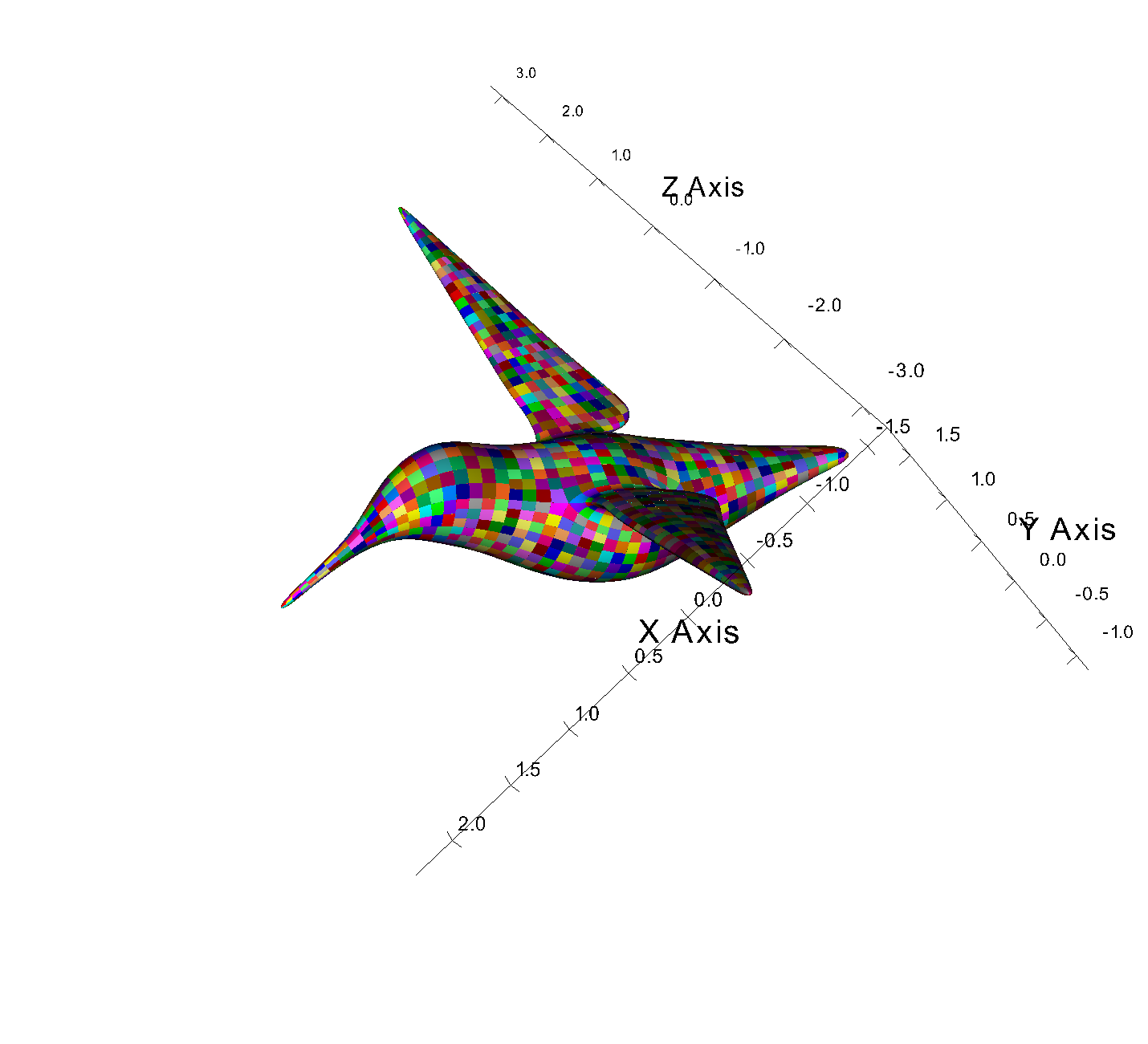}
\caption{\footnotesize Patch configuration for the hummingbird model in the setup of
numerical experiment described in Section~\ref{sec:exm_hummingbird}.\label{fig:hummingbird_patch}}
\end{figure}

\begin{figure}[ht]
  \centering
\includegraphics[scale=0.26, trim={240 100 70 220}, clip]{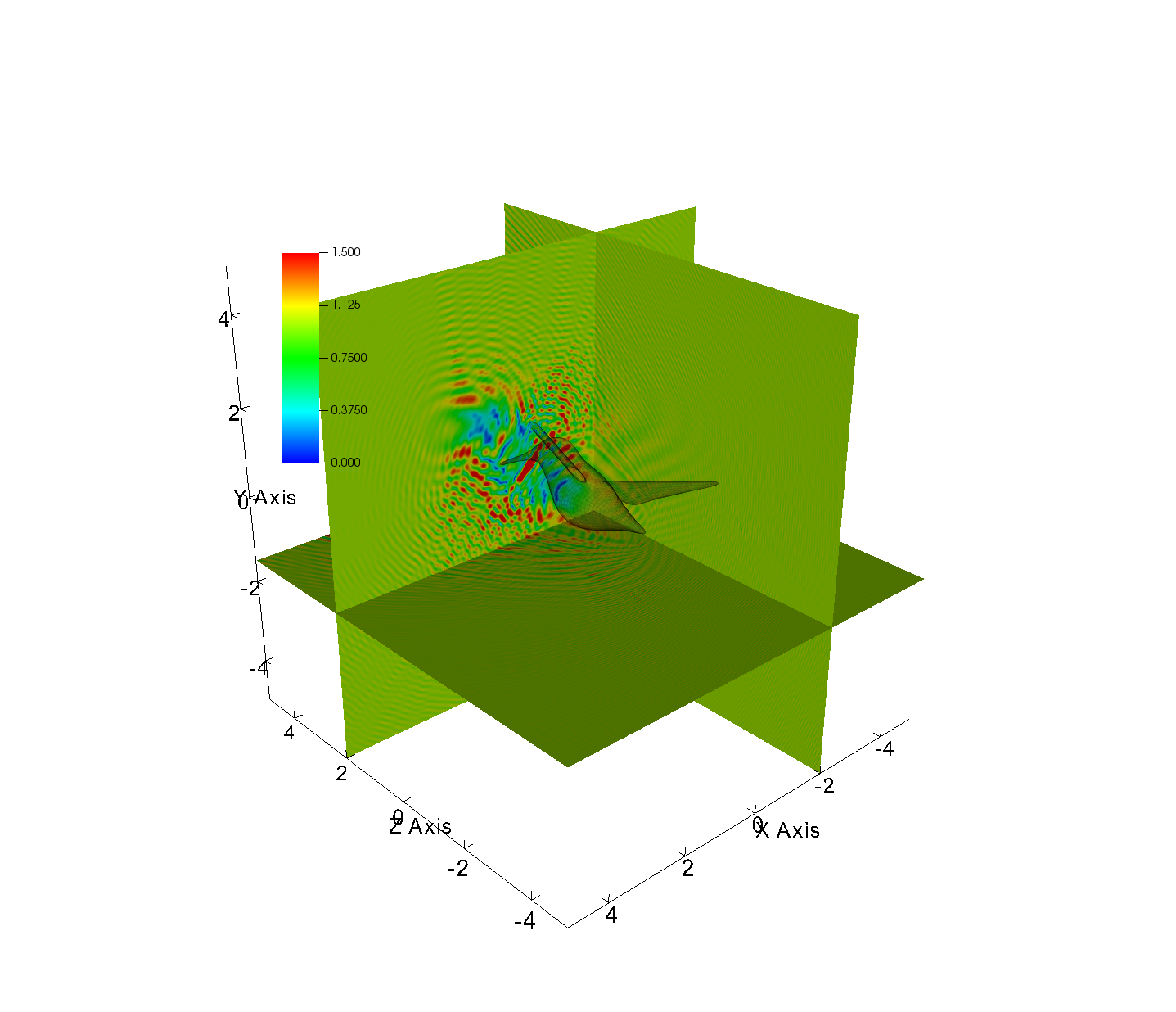}
\caption{\footnotesize Absolute value of the field component \(E_{x}\) on the
cross-section of the \([-5,5]\times [-5,5]\) squares on the planes
\(x=-2\), \(y=-1.5\), and \(z=2\). Field values are scaled
to \([0,1.5]\).\label{fig:hummingbird_scaled_abs}}
\end{figure}

Next, we consider the NURBS hummingbird model
\cite{cad_hummingbird_src} with
the following material properties:
\(\epsilon_{e}=\mu_{e}=1,\kappa_{e}=50.27\), and 
\(\epsilon_{i}=2.16,\mu_{i}=1,\kappa_{i}=73.89\) with the wingspan of
the scatterer being \(77\) wavelengths. The refractive index of the medium
in \(\Omega_i\) is \(1.47\). The discretization contains
\(N\approx 680k\) points with \(P=4728\) curvilinear patches and
\(N_{C}=12\) points in each variable. The patch decomposition of the surface of the hummingbird model for this simulation is shown in Fig.~\ref{fig:hummingbird_patch}. The accelerated solver takes
\(27\) seconds on average in one iteration against \(588\) seconds (an
average of five runs) by the unaccelerated solver.
We plot the absolute value of the field component \(E_{x}\) of the computed
total field on the union of the truncated 
planes \(x=-2\), \(y=-1.5\), and \(z=2\) in
Fig.~\ref{fig:hummingbird_scaled_abs}.
\subsection{Scattering From Nanophotonic Splitter CAD Model}\label{sec:exm_splitter}
\begin{figure}[ht!]
\includegraphics[scale=0.19]{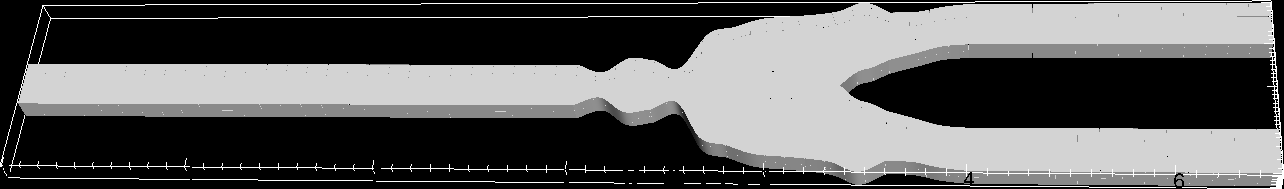}
\caption{Shape of the nanophotonic splitter scatterer.\label{fig:splitter}}
\end{figure}

\begin{figure}[ht!]
  \centering
  \includegraphics[scale=0.2]{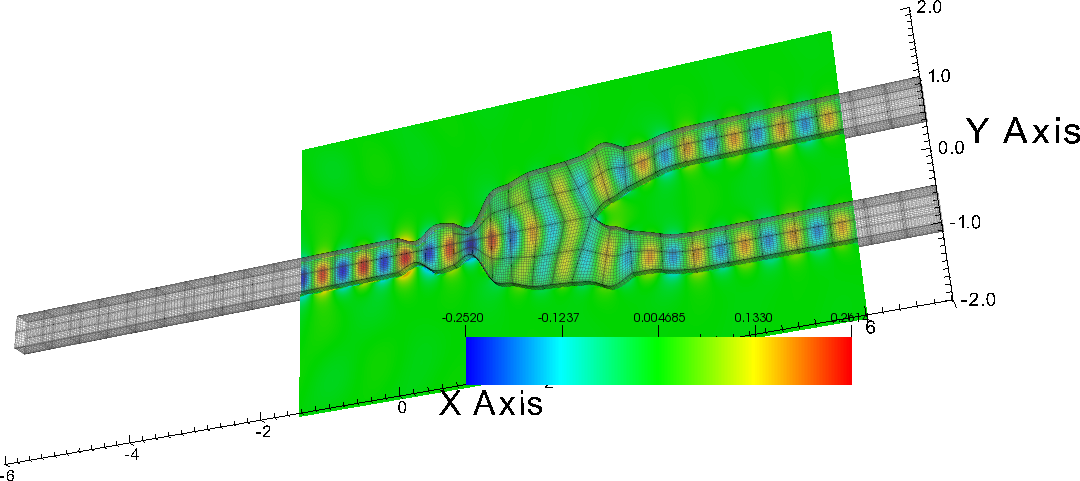}
\caption{\footnotesize Real part of \(E_{y}\).\label{fig:splitter_ey_re}}
\end{figure}

\begin{figure}
\centering
  \includegraphics[scale=0.2]{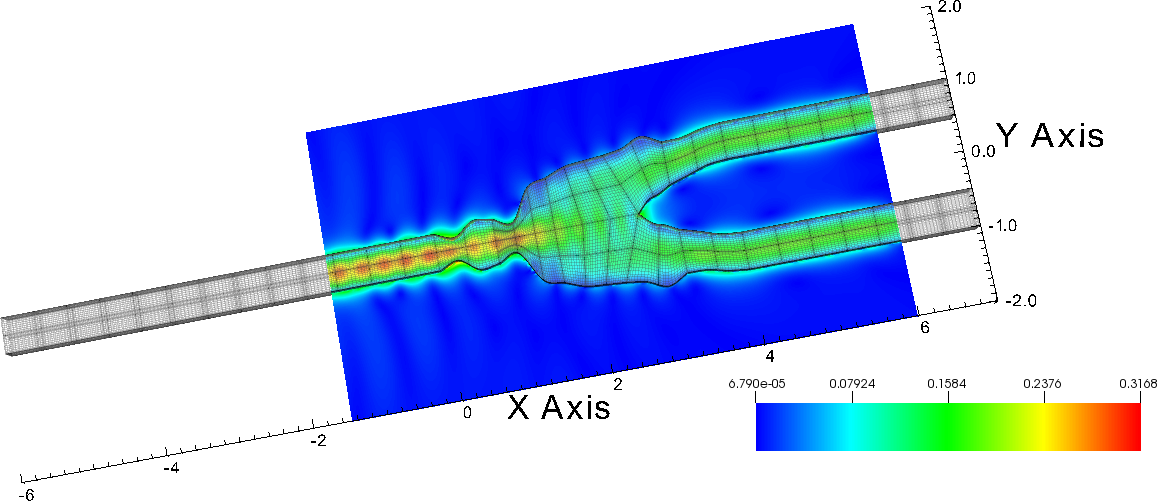}
\caption{\footnotesize Absolute value of \(E_{y}\).\label{fig:splitter_ey_abs}}
\end{figure}

For the final numerical example, we consider a
\emph{Gmsh}-produced~\cite{geuzaine_gmsh} geometry of a nanophotonic
splitter~\cite{garza_acs_2023} (with size \(29\lambda_{i}\) in length) as depicted in
Fig.~\ref{fig:splitter}. The material properties considered for
this simulation are as follows: 
(\(\epsilon_{e}=2.08\), \(\mu_{e}=1\), \(\kappa_{e}=5.85\)) and
(\(\epsilon_{i}=12.08\), \(\mu_{i}=1\), \(\kappa_{i}=14.09\)). The
refractive indices within the domains \(\Omega_{e}\) and
\(\Omega_{i}\) are \(1.444\) and \(3.476\), representing silicon dioxide and silicon respectively. An electric dipole was used to excite the input waveguide of the splitter.  
The surface is discretized with \(P=280\) curvilinear patches with
\(N_{C}=12\) points per patch in each variable within each patch containing a total of
\(N\approx 40k\) points.
For simulation in this experiment, in conjunction with the proposed
method, we have additionally used a windowing function
to properly truncate the integral operators, see~\cite{garza_acs_2023} and the relevant
references therein. The run time for one FM for the unaccelerated CBIE
and the CBIE-IFGF methods are \(2.49\) seconds and \(0.52\) seconds,
respectively, demonstrating that the CBIE-IFGF provides a 5x speed-up even for this relatively small-size problem. Fig.~\ref{fig:splitter_ey_re} and
Fig.~\ref{fig:splitter_ey_abs} show the real part and the absolute
value, respectively, of the component \(E_{y}\) of the computed total field on the rectangle
\((x,y)\in[-1.5,6]\times[-2,2]\) on the \(xy\)-plane.
\section{Summary}\label{sec:summary}
This paper presents an accelerated Nystr\"{o}m solver for the
electromagnetic scattering problem for a dielectric medium. The method
uses non-overlapping curvilinear parametric patches to decompose the surface, and
approximates the density on a tensor-product of 1D Chebyshev
grids in both the parametric variables. A polynomial change-of-variables-based strategy is used to resolve the singular behavior of the Green's function in conjunction with the IFGF interpolation strategy to
achieve an overall \(\mc{O}(N\log N)\) computational complexity. We
illustrated the accuracy and the efficiency of the proposed method with multiple numerical examples. Future work includes further acceleration of the EM
scattering problem using GPU programming.

\end{document}